\input amstex
\loadbold
\loadeurm
\loadeusm
\openup.8\jot \magnification=1200

\def\L{\Lambda}
\def\O{\Omega}
\def\cls{{\text{ cls}}}
\def\Z{\Bbb Z}
\def\A{\Bbb A}
\def\Q{\Bbb Q}

\def\t{\theta}

\def\rank{{\text {rank}}}
\def\mod{{\text{mod\ }}}
\def\mass{{\text{mass}}}
\def\D{\Delta}
\def\H{{\Bbb H}}
\def\rad{{\text {rad}}}
\def\e{{\text{e}}}

\def\R{{\Bbb R}}
\def\C{{\Bbb C}}
\def\Tr{{\text{Tr}}}
\def\mult{{\text{mult}}}

\def\gen{{\text{gen}}}
\def\cls{{\text{cls}}}
\def\F{{\Bbb F}}
\def\char{{\text{char}}}

\documentstyle{amsppt}
\pageheight{7.7in}
\vcorrection{-0.05in}
\topmatter
\pageno 1
\title Action of Hecke Operators on Siegel Theta Series II\endtitle
\author Lynne H. Walling\endauthor
\address L.H. Walling, Department of Mathematics, University of Bristol,
Bristol BS8 1TW, England
\endaddress
\email l.walling\@bristol.ac.uk \endemail

\abstract We apply the Hecke operators $T(p)^2$ and  $T'_j(p^2)$
($1\le j\le n\le 2k$) to a degree $n$ theta series attached to a 
rank $2k$ $\Z$-lattice $L$
equipped with a positive
definite quadratic form in the case that $L/pL$ is
regular.  We 
explicitly realize the image of the theta series under these Hecke
operators as a sum of theta series attached to certain 
sublattices of ${1\over p}L$,
thereby generalizing the Eichler Commutation Relation.
We then show that the average
theta series (averaging over isometry classes in a given genus)
is an eigenform for these operators.  
We explicitly compute the eigenvalues on the average theta series,
extending previous work where we had the restrictions that $\chi(p)=1$
and $n\le k$.  We also show that $\t(L)|T'_j(p^2)=0$ for
$j>k$ when $\chi(p)=1$, and
for $j\ge k$ when $\chi(p)=-1$, and that $\t(\gen L)$ is an eigenform for
$T(p)^2$.  
\endabstract

\endtopmatter
\document

\head{\S1. Introduction and statements of results}\endhead

The Fourier coefficients of a degree $n$ Siegel theta series tell us how
many times a given positive definite quadratic form of
rank $2k$ over
$\Z$ represents each rank $n$ quadratic form.  Hecke operators help us
study Fourier coefficients of modular forms.

In this paper we complete the analysis begun in [12], examing
the action of the Hecke operators on the Fourier coefficients of a
Siegel theta series of degree $n$.  
We first extend the Eichler Commutation Relation,
describing the image of the theta series $\t(L)$ under the (below defined)
Hecke operators $T'_j(p^2)$ ($1\le j\le n\le 2k$) as a sum of theta series attached
to certain sublattices of ${1\over p}L$.  Then averaging over the genus
of $L$ (see the definition later in this section), we find that 
for $j\le k$, $\t(\gen L)$ is an eigenform for
$T'_j(p^2)$, with eigenvalue
$$\lambda_j(p^2)=\cases
p^{j(k-n)+j(j-1)/2} \beta(n,j) (p^{k-1}+1)\cdots(p^{k-j}+1)
&\text{if $\chi(p)=1$,}\\
p^{j(k-n)+j(j-1)/2} \beta(n,j) (p^{k-1}-1)\cdots(p^{k-j}-1)
&\text{if $\chi(p)=-1$,}
\endcases$$
where $\beta(n,j)$ is the number of $j$-dimensional subspaces of an
$n$-dimensional space over $\Z/p\Z$.

We also show that $\t(L)$ vanishes under $T'_j(p^2)$ for $j>k$ when
$\chi(p)=1$, and for $j\ge k$ when $\chi(p)=-1$.  Using this, we show
$\t(\gen L)$ is an eigenform for $T(p)^2$ with eigenvalue
$$\lambda(p^2) 
=\cases \left((p^{k-1}+1)\cdots(p^{k-n}+1)\right)^2 &\text{if $\chi(p)=1$,}\\
\left((p^{k-1}-1)\cdots(p^{k-n}-1)\right)^2 &\text{if $\chi(p)=-1$.}\endcases$$

Let us now state our assumptions, present the relevant definitions, and
outline our strategy.

Throughout, $L$ is a rank $2k$ lattice over $\Z$ equipped with
a positive definite quadratic form $Q$.  By scaling $Q$ if necessary,
we can assume $L$ is even integral, meaning
$Q(L)\subseteq 2\Z$.  With $B$ the symmetric bilinear form
associated to $Q$ so that $Q(v)=B(v,v)$, $(v_1,\ldots,v_{2k})$ a
$\Z$-basis for $L$, and $A=\big( B(v_i,v_j)\big)$, we have
$$Q(\alpha_1 v_1+\cdots+\alpha_{2k}v_{2k})
= (\alpha_1\cdots\alpha_{2k})A\ ^t(\alpha_1\cdots\alpha_{2k}).$$

The quotient $L/pL$ is a vector space over $\Z/p\Z$, with induced quadratic
form $Q$ modulo $p$ when $p$ is odd, $Q'={1\over 2}Q$ modulo 2 when $p=2$.
A subspace $\overline C$ of $L/pL$ is called totally isotropic if all
its vectors vanish under the induced quadratic form.

Given another lattice $K$ on the space $\Q L$, we use $\{L:K\}$ to denote
the invariant
 factors, also called the elementary divisors, of $K$ in $L$
(see \S81D of [8]).  We write $\mult_{\{L:K\}}(a)$ to denote
the multiplicity of $a$ as an invariant factor of $K$ in $L$.

The Siegel theta series attached to $L$ is
$$\t(L;\tau) = \sum_C \e\{\ ^tCAC\tau\}$$
where $C$ varies over $\Z^{2k,n}$, 
$\tau\in \{X+iY:\ \text{symmetric } X,Y\in\R^{n,n},\ Y>0\ \},$
and $\e\{*\} = \exp(\pi i \Tr(*)).$
(Here $Y>0$ means that the quadratic form represented by the matrix $Y$
is positive definite.)
Since $Q$ is positive definite and $\Im \tau>0$, the series $\t(L;\tau)$
is absolutely convergent.  
We also set $\t(\gen L)=\sum_{L'} {1\over o(L')} \t(L')$ where
$L'$ varies over the isometry classes in the genus of $L$, and $o(L')$
is the order of the orthogonal group of $L'$.  ($L'$ is in the genus
of $L$ if, locally everywhere, $L'$ and $L$ are isometric.)
Note that some authors normalize
this average to have 0-coefficient equal to 1.

As $C$ varies, $(v_1,\ldots,v_{2k})C$ varies over all $(x_1,\ldots,x_n)$,
$x_i\in L$.  Let $\L$ be the (formal) direct sum
$\Z x_1\oplus\cdots\oplus\Z x_n$ equipped with the (possibly semi-definite)
quadratic form given by $T=(B(x_i,x_j))$.  Let
$$\e\{\L\tau\}=\sum_G \e\{\ ^tGTG\tau\}$$
where $G$ varies over $GL_n(\Z)$ (or, if $k$ is odd, we equip $\L$ with
an orientation and let $G$ vary over $SL_n(\Z)$).
Then as the $\L$ vary over all formally rank $n$ sublattices of $L$, we have
$$\t(L;\tau) = \sum_{\L} \e\{\L\tau\}.$$

Note that $\t(L;\tau) = \sum_T r(A,T) \e\{ T\tau\}$ where
$$r(A,T) = \#\{C\in\Z^{2k,n}:\ ^tCAC=T\ \}.$$
Here $\rank \ ^tCAC$ can never exceed $2k$, which is why
 we restrict our attention
to $n\le 2k$.  (See chapter IV of [6] to read about ``singular''
Siegel modular forms, which are series whose support contain only
singular matrices.)

Also, $\t(L)$ ``transforms'' under a congruence
subgroup of 
$$\align &Sp_n(\Z)\\
&\qquad = \left\{\pmatrix A&B\\C&D\endpmatrix \in GL_{2n}(\Z):\ 
A\ ^tB,\ C\ ^tD \text{ symmetric},\ A\ ^tD - B\ ^tC = I \ \right\}.
\endalign$$
More precisely,
$$\t\big(L;(A\tau+B)(C\tau+D)^{-1}\big)
= \chi(\det D)\det(C\tau+D)^k\t(L;\tau)$$
for all $\pmatrix A&B\\C&D\endpmatrix \in Sp_n(\Z)$ with $C\equiv 0\ (\mod N).$
Here $N$ is the ``level'' of $L$, i.e. the smallest positive integer so that
$NA^{-1}$ is even integral (meaning $NA^{-1}$ is an integral matrix
with even diagonal).  Also, $\chi$ is a quadratic
 Dirichlet character modulo $N$.
In fact, given a prime $p$, $p|N$ if and only if $L/pL$ is ``not regular''
(meaning $L/pL$ has a nontrivial 
totally isotropic subspace orthogonal to all of $L/pL$).  
Also, for $p\nmid N$, $\chi(p)=1$ if and only if $L/pL$ is 
``hyperbolic'', and $\chi(p)=-1$ otherwise.  (A hyperbolic plane is a
dimension 2 space with quadratic form given by the matrix 
$\pmatrix 0&1\\1&0\endpmatrix$; a space is hyperbolic if it is the
orthogonal sum of hyperbolic planes.)

For each prime $p$,  associated to $p$
 there are $n+1$ Hecke operators $T(p)$,
$T_j(p^2)$ ($1\le j\le n$).
In [5], we analyzed the action on Fourier coefficients of the operators
$T(p)$ and $\widetilde T_j(p^2)$
(where the $\widetilde T_j(p^2)$ are simple linear combinations of
$T_{\ell}(p^2)$, $0\le \ell\le j$; see Theorem 1.1 (b) below).  
We did this by finding a set of coset
representatives for these operators.  Then Theorem 6.1 of [5] states:

\proclaim{Theorem 1.1} Let $F$ be a Siegel modular form of degree $n$,
weight $k$, level $N$, and character $\chi$, and expand $F$ as
$$F(\tau) = \sum_{\L} c(\L)\e^*\{\L\tau\}$$
where $\L$ varies over even integral positive semi-definite isometry
classes of rank $n$ lattices, and $\e^*\{\L\tau\}=\sum_G\e\{\ ^tGTG\tau\}$
where $T$ is a lattice giving the quadratic form on $\L$, and $G$ varies
over $O(T)\backslash GL_n(\Z)$ when $k$ is even,
$O^+(T)\backslash SL_n(\Z)$ when $k$ is odd.  (Here $O(T)$ denotes
the orthogonal group of $T$.  Also, when $k$ is odd, we equip $\L$ with
an orientation.)
\item{(a)}  The coefficient of $\e^*\{\L\tau\}$ in $F|T(p)$ is
$$\sum_{p\L\subseteq\O\L} \chi([\O:p\L]) p^{E(\L,\O)}c(\O^{1/p})$$
where $E(\L,\O)=m(1)k+m(p)(m(p)+1)/2-n(n+1)/2$, $m(a)=\mult_{\{\L:\O\}}(a)$,
and $\O^{1/p}$ denotes that lattice $\O$ scaled by $1/p$.
\item{(b)}  For $1\le j\le n$, set 
$$\widetilde T_j(p^2)=p^{j(k-n-1)}\sum_{0\le\ell\le j} \beta(n-\ell,j-\ell)
T_j(p^2)$$ where $\beta(m,r)=\prod_{i=0}^{r-1} {p^{m-i}-1\over p^{r-i}-1}$,
(so this is
the number of $r$-dimensional subspaces of an $m$-dimensional space over
$\Z/p\Z$ when $m\ge r$).
Then the coefficient of $\e^*\{\L\tau\}$ in
$F|\widetilde T_j(p^2)$ is
$$\sum_{p\L\subseteq\O\subseteq{1\over p}\L}
 \chi(p^{j-n}[\O:p\L]) p^{E_j(\L,\O)} \alpha_j(\L,\O) c(\O)$$
where $E_j(\L,\O)= k(m(1/p)-m(p)+j)+m(p)(m(p)+m(1)+1)
+m_j(1)(m_j(1)+1)/2 - j(n+1) $,
$m_j(1)=m(1)-n+j$,
and $\alpha_j(\L,\O)$ denotes the number of totally isotropic co-dimension
$n-j$ subspaces of $(\L\cap\O)/p(\L+\O)$.
\endproclaim

\noindent{\bf Remark.}
Above we wrote a Siegel modular form
as a series in terms of $\e^*\{\L\tau\}$, whereas
we previously wrote $\t(L)$ as a series in $\e\{\L\tau\}$,
$\L=\Z x_1+\cdots+\Z x_n\subseteq L$.  Letting $O(\L)$ denote the orthogonal
group of $\L$ as a positive-definite sublattice of $L$ with
$\rank\L\le n$, 
and $o(\L)=\#O(\L)$, one can show $\e\{\L\tau\}=o(\L)\e^*\{\L\tau\}$.
\smallskip

The strategy used here is essentially the same as that used in [12]
where we were restricted to $\chi(p)=1$ and $j\le n\le k$.
As in Proposition 1.4 of [12], in Proposition 2.1 we first directly
apply to $\t(L;\tau)$ the matrices found in Corollary
2.1 of [5] that give the action of the operators $\widetilde T_j(p^2)$
for $\chi(p)=\pm1$, $n\le 2k$.  We find that $\widetilde c_j(\O)$,
the  coefficient of $\e\{\O\tau\}$ in $\t(L)|\widetilde T_j(p^2)$, 
is a sum over $\L$ where $p\O\subseteq\L
\subseteq({1\over p}\O\cap L)$; we construct all these $\L$ while
simultaneously computing the summand attached to $\L$.  (In [12],
$\widetilde c_j(\O)$ was called $c^*_j(\O)$; we hope this revised
notation is more suggestive of the quantity represented.)
Then as in Proposition 1.5 of [12], in Proposition 2.2 we compute
$b_j(\O)$, the coefficient of $\e\{\O\tau\}$ in $\sum_{K_j}\t(K_j)$ where
$K_j$ varies over all lattices in $\gen L$ with $pL\subseteq K_j\subseteq
{1\over p}L$, $\mult_{\{L:K_j\}}(1/p)=\mult_{\{L:K_j\}}(p)=j$,
$\chi(p)=\pm1$ and $n\le 2k$; the geometry of $L$ constrains this
computation to $j\le k$ when $\chi(p)=1$, $j<k$ when $\chi(p)=-1$.
Then as in Theorem 1.2 of [12], in Theorem 2.3 we 
use these propositions to realize $\t(L)|T'_j(p^2)$
as a linear combination of $\t(K_{\ell})$, where $0\le \ell\le j$ and
$T'_j(p^2)$ is a specific linear combination of the $\widetilde T_{\ell}(p^2)$,
$0\le\ell\le j$ (defined in Theorem 2.3); here $j\le k$ if $\chi(p)=1$,
$j<k$ if $\chi(p)=-1$.

In Corollary 2.4 we extend
Corollary 1.3 of [12], showing the average theta
series $\t(\gen L)$ is an eigenform for the $T'_j(p^2)$
where $j\le k$ when $\chi(p)=1$, $j<k$ when $\chi(p)=-1$; 
we explicitly
compute the eigenvalues.

In Proposition 3.1, we consider $j>k$ when $\chi(p)=1$, $j\ge k$ when
$\chi(p)=-1$.  We realize $\t(L)|\widetilde T_j(p^2)$ as a linear combination
of $\t(L)|\widetilde T_{\ell}(p^2)$, $\ell\le k$ when $\chi(p)=1$,
$\ell<k$ when $\chi(p)=-1$.  Then in Theorem 3.3, we show
$\t(L)|T'_j(p^2)=0$ for $j>k$ when $\chi(p)=1$, $j\ge k$ when
$\chi(p)=-1$.
Finally, in Theorem 3.4 we use the preceeding results
and the formula from Proposition 5.1 of [5] realizing $T(p)^2$ as a linear
combination of $\widetilde T_j(p^2)$, $0\le j\le n$, to show $\t(\gen L)$
is an eigenform for $T(p)^2$, explicitly computing the eigenvalue.

In \S4 we extend Lemma 1.6 of [12] and collect some useful combinatorial 
identities.

In many of our arguments we work in a quadratic space over a finite field
$\F$ with characteristic $p$.  When $p\not=2$, we directly apply theorems
from \S42 and \S62 of [8].
When $p=2$, we could use the 
results on lattices over local dyadic rings in \S93
of [8] to deduce the results we need;
for completeness, in \S5 we give a self-contained treatment of quadratic
spaces over finite fields with characteristic 2.

The proofs of Propositions 2.1 and 2.2 closely parallel those of Propositions
1.4 and 1.5 of [12].  The main technique is to construct and count
lattices $\L$, $p\O\subseteq\L\subseteq{1\over p}\O$ (where $\O$ is
given), so that we control the structure of $\L$.  We do this by using a
two-step modulo $p$ construction; the constructions differ in these two
propositions, but the approach is the same.

To construct and count the $K_j$ of Proposition 2.2, we first construct
a dimension $j$ totally isotropic subspace $\overline C$ of $L/pL$
(which is a vector space over $\Z/p\Z$ with quadratic form $Q$ modulo $p$
if $p$ is odd, and quadratic form ${1\over 2}Q$ modulo $2$ when $p=2$).
We set $K'$ equal to the preimage of $\overline C$ in $L$.  Knowing
the structure of $\Z_p L$, we infer the stucture of $\Z_p K'$.
Then in $K'/pK'$ (scaled by $1/p$), we refine $\overline C$, building
$\overline C'$, a dimension $j$ totally isotropic subspace independent
of $\overline{pL}$; we set $pK_j$ equal to the preimage in $K'$ of
$(\overline C')^{\perp}$, the orthogonal complement of $\overline C'$.
Thus we control both the local structure of $K_j$ and its invariant
factors in $L$.

To construct and count the rank $n$ lattices $\L\subseteq L$
of Proposition 2.1 where
$\O\subseteq{1\over p}L$ is fixed with rank $n$ and
$p\O\subseteq \L\subseteq{1\over p}\O$, we first note that we must have
$\L=\O_0\oplus\L'$, $\L'\subseteq\O_1\oplus\O_2$, where
$\O={1\over p}\O_0\oplus\O_1\oplus p\O_2$, $\O_i\subseteq L$ with
$\O_1\oplus\O_1$ primitive modulo $p$ in $L$ (meaning 
$(\O_0\oplus\O_1)\cap pL=p(\O_0\oplus\O_1)$).  So in this two-step process,
we begin with
$$\D={1\over p}\O\cap L=\O_0\oplus\O_1\oplus\O_2,$$
then in $\D/p\D$ we extend $\overline{\O\cap\D}=\overline{\O_0\oplus\O_1}$
to $\overline{\O_0\oplus\O_1\oplus\D_2}$ where we control 
$\dim \overline\D_2$. Letting $\D'$ be the preimage in $\D$ of
$\overline{\O_0\oplus\O_1\oplus\D_2}$, in $\D'/p\D'$ we extend
$\overline{p\O}=\overline \O_0$ to $\overline{\O_0\oplus U}$ where
$\overline U$ is totally isotropic of a given
dimension $\ell$; this will enable us to simultaneously compute
$\alpha_j(\L,\O)$ as we construct $\L$.  Then we extend $\overline{\O_0
\oplus U}$ to $\overline{\O_0\oplus\L_1\oplus\L_2}$ where $\overline U
\subseteq\overline\L_1$, $\overline\L_1$ is independent of $\overline{p\D}$,
$\overline\L_2$ is independent of $\overline{\O\cap\D}$, and we specify
$\dim\overline \L_i$.  Consequently, letting $\L$ be the preimage
in $\D'$ of $\overline{\O_0\oplus\L_1\oplus\L_2}$, we get
$$\L=\O_0\oplus(\L_1\oplus p\L_1')\oplus (\L_2\oplus p\L_2'\oplus p^2\L_2'')$$
where $\O_1=\L_1\oplus\L_1'$, $\O_2=\L_2\oplus\L_2'\oplus\L_2''$.
The quantity $\alpha_j(\L,\O)$ counts the number of totally isotropic
codimension $n-j$ subspaces of $(\L\cap\O)/p(\L+\O)\approx
\L_1/p\L_1\oplus p\L_2'/p^2\L_2'$.
Thus the subspaces counted by $\alpha_j(\L,\O)$ that project onto a given
$\overline U$ of dimension $\ell$ in $\L_1/p\L_1$ is the number of dimension
$d-\ell$ subspaces of $p\L_2'/p^2\L_2'$ where $d=d_1+d_2'-n+j$,
$d_1=\rank\L_1$, $d_2'=\rank\L_2'$.
Hence our construction allows us to
 control the invariant factors of $\L$ in $\O$ as we compute $\alpha_j(\L,\O)$.

The proof of Theorem 2.3 relies on Lemma 4.1, an extension of the
Reduction Lemma (Lemma 1.6) of [12].  Lemma 4.1 allows us to write
$\varphi_j(U\perp\H^t)$ and $\varphi_j(U\perp\H^t\perp\A)$ in terms
of $\varphi_{\ell}(U)$, $\ell\le j$; here $U$ is a quadratic space over
$\Z/p\Z$, $\H\simeq\pmatrix 0&1\\1&0\endpmatrix$ is a hyperbolic plane,
$\A$ is an anisotropic plane, and $\varphi_{\ell}(U)$ denotes the number
of $\ell$-dimensional totally isotropic subspaces of $U$.
(Recall that a space $W$ is totally isotropic if $Q(W)=0$.)  In [12],
Lemma 1.6 relates $\varphi_j(U\perp\H^t)$ to $\varphi_{\ell}(U)$,
$\ell\le j$, when $t$ is positive.  In Lemma 4.1, we handle
``cancellation'' of an anisotropic plane, and $t$ negative (when this
is meaningful; see discussion preceeding Lemma 4.1).  
Also, when $\chi(p)=-1$, 
we could not formulate $b_j(\O)$ in 
the same way we did when $\chi(p)=1$ (compare
Proposition 2.2 herein with Proposition 1.5 of [12]). Thus the 
argument used to prove Theorem 2.3 is not identical to that used to
prove Theorem 1.2 of [12], although both proofs are simple applications of
basic combinatorial identities.

The reader is referred to [1], [3], and [7] for facts about Siegel modular
forms, and to [2] and [8] for facts about quadratic forms.
See, for instance, [1], [4], [9], [10], [11], [12], [13], [14] for earlier
work on the (generalized) Eichler Commutation Relation and the action
of Hecke operators on theta series.

\head{\S2. Nonzero eigenvalues of Hecke operators $T'_j(p^2)$}
\endhead
Throughout, $L$ is an even
integral, rank $2k$ lattice with positive definite quadratic
form of level $N$; we
fix a prime $p$, $p\nmid N$.
For $0< r$ and any $m$,  set
$$\delta(m,r)=\epsilon(m-r+1,r)=\prod_{i=0}^{r-1}(p^{m-i}+1),\ 
\mu(m,r)=\prod_{i=0}^{r-1}(p^{m-i}-1),$$
and
$$\beta(m,r)=\prod_{i=0}^{r-1}{p^{m-i}-1\over p^{r-i}-1}
={\mu(m,r)\over \mu(r,r)}.$$
We agree that when $r=0$, the value of any of these functions is 1.
(In [12] we used the function $\epsilon$; here we use instead the function
$\delta$
 which we define as a product indexed as is the product defining $\beta$,
allowing us to more readily see similarities between $\delta$ and $\beta$.)

Let $V$ be a quadratic space over $\Z/p\Z$ with quadratic form $Q$
(e.g. $V=L/pL$).  
The radical of $V$ is
$$\rad V=\{x\in V:\ Q(x)=0 \text{ and }B(x,V)=0\ \}$$
where $B$ is the symmetric bilinear form associated to $Q$ so that
$Q(x)=B(x,x)$ if $p\not=2$, and $Q(x)={1\over 2}B(x,x)$ if $p=2$.
It is easily seen that if
$V=U\perp\rad V=U'\perp\rad V$ then $U$ is isometric to $U'$, written
$U\simeq U'$.
We say $V$ is regular if  $\rad V=\{0\}.$
A subspace $U$ of $V$ is called a hyperbolic
plane if it has dimension 2 and its quadratic form is given by the matrix
$\pmatrix 0&1\\1&0\endpmatrix$; $U$ is called hyperbolic if it is the 
orthogonal sum
of hyperbolic planes. A nonzero vector $x\in V$ is called isotropic if
$Q(x)=0$.  A space is called isotropic if it contains at least one 
(nonzero) isotropic
vector, and anisotropic otherwise; a space is called totally isotropic
if all its (nonzero) vectors are isotropic.
Also, by 62:1a of [8], a regular space over $\Z/p\Z$, $p$ odd, is completely
determined by its dimension and the square class of its discriminant.
For an analogous result when $p=2$, see the discussion following
Proposition 5.4.

We let $\varphi_{\ell}(V)$ denote the
number of $\ell$-dimensional totally isotropic subspaces of $V$.
 When $p$ is odd,
we rely on formulas from p. 143-146 [2] that give us
$\varphi_1(U)$ when $U$ is regular; when $p=2$,
we use Theorem 5.11.  These formulas show (see [11]) that
when $U$ is regular,
$$\varphi_{\ell}(U)=
\cases
\beta(t,\ell)\delta(t-1,\ell)
&\text{if $\dim U=2t$ and $U$ is hyperbolic,}\\
\beta(t-1,\ell)\delta(t,\ell)
&\text{if $\dim U=2t$ and $U$ is not hyperbolic,}\\
\beta(t,\ell)\delta(t,\ell)
&\text{if $\dim U=2t+1$.}
\endcases$$

We have 
$$\t(L;\tau)=\sum_{\Lambda} \e\{\Lambda \tau\}$$
where $\Lambda$ varies over all sublattices of $L$ with (formal) rank $n$.
(So $\Lambda$ is the external direct sum $\Z x_1\oplus\cdots\oplus\Z x_n$
where $x_1,\ldots,x_n \in L$.)

Let $p$ be a prime not dividing $N$, the level
of $L$.  We have $n+1$ Hecke operators associated to $p$, named $T(p)$,
$T_j(p^2)$ ($1\le j\le n$).  For $T$ one of these operators, there is an
associated matrix $\delta$ so that
$$F|T=p^{\eta} \sum_{\gamma} F|\delta^{-1}\gamma$$
where $\gamma$ runs over $(\Gamma\cap\Gamma')\backslash\Gamma$,
$\Gamma=\Gamma_1(N)$, $\Gamma'=\delta\Gamma\delta^{-1}$, and
$p^\eta$ is a normalizing factor.  Here $\delta=\pmatrix pI_n\\&I_n\endpmatrix$
and $\eta=n(k-n-1)/2$ when $T=T(p)$; 
$$\delta=\pmatrix pI_j\\&I_{n-j}\\
&&{1\over p}I_j\\&&&I_{n-j}\endpmatrix$$
 and $\eta=0$ when
$T=T_j(p^2)$.
As discussed in [5], when analyzing the action of the operators $T_j(p^2)$
on Fourier coefficients of Siegel modular forms, we encounter incomplete
character sums.  To complete these character sums, we set
$$\widetilde T_j(p^2)=p^{j(k-n-1)}\sum_{0\le\ell\le j} \beta(n-\ell,j-\ell)
T_j(p^2)$$ where $\beta(m,r)=\prod_{i=0}^{r-1} {p^{m-i}-1\over p^{r-i}-1}$.

In Propositions 2.1 and 3.1 of [5] we find explicit coset representatives
giving the action of each of these operators.
We directly apply the coset representatives for $\widetilde T_j(p^2)$
to get the following.

\proclaim{Proposition 2.1}
Let $1\le n\le 2k$, $1\le j\le n$, and $p$ a prime so that $p\nmid N$
($N$ the level of $L$).
Write
$$\t(L;\tau)|\widetilde T_j(p^2) = \sum_{\O} \widetilde c_j(\O) \e\{\O\tau\}$$
where $\O$ varies over even integral sublattices of ${1\over p}L$ that
have (formal) rank $n$.

\item{(a)}  Say $\chi(p)=1$.  Then
$$\widetilde c_j(\O)
 = \sum_{\ell,t}p^E \varphi_{\ell}(\overline\O_1)\delta(k-r_0-\ell-1,t)
\beta(r_2,t)\beta(n-r_0-\ell-t,j-r_0-\ell-t)$$
where $E=E'(\ell,t,\O)
=\ell(k-r_0-r_1)+\ell(\ell-1)/2+t(k-n)+t(t+1)/2$.
\item{(b)}  Say $\chi(p)=-1$.  Then
$$\widetilde c_j(\O)
 = \sum_{\ell,t}(-1)^{\ell}
p^E \varphi_{\ell}(\overline\O_1)\beta(k-r_0-\ell-1,t)
\mu(r_2,t)\beta(n-r_0-\ell-t,j-r_0-\ell-t)$$
where $E=E'(\ell,t,\O)$ is as in (a).
\endproclaim

\demo{Proof}  (a) In Proposition 1.4 of [12] we proved a formula for
$\widetilde c_j(\O)$ (there called $c^*_j(\O)$) provided $\chi(p)=1$.
Making the change of variables $t\mapsto j-r_0-t$ yields (a).

(b)
As in the proof of Proposition 1.4 of [12], we directly apply to $\t(L)$
coset representatives giving the action of $\widetilde T_j(p^2)$.  The
coset representatives we use are from Proposition 2.1 of [5]
(see also Corollary 2.1 and Theorem 4.1 of [5]).  In Theorem 6.1 of
[5] we used these coset representatives to
examine the action of Hecke operators on Siegel modular forms with level
and character.  
So, applying our coset representatives to $\t(L)$ we initially get
$$\align
\t(L;\tau)|\widetilde T_j(p^2)
&=\sum_{\L\subseteq L}\left( \sum_{{p\L\subseteq\O\subseteq{1\over p}\L}\atop
{\O\text{\ integral}}} 
p^{E_j(\O,\L)} \alpha_j(\L,\O) \right) e\{\O\tau\} \\
&= \sum_{{\O\subseteq {1\over p}L}\atop{\O\text{\ integral}}}
\left( \sum_{p\O\subseteq\L\subseteq({1\over p}\O\cap L)} 
(-1)^{j-m_0+m_2}
p^{E_j(\O,\L)} \alpha_j(\L,\O)\right) e\{\O\tau\} \endalign$$
where
$$\align
E_j(\L,\O)&
=k(j-m_0+m_2) +m_0(n-m_2+1)\\
&\qquad  +(j-m_0-m_2)(j-m_0-m_2+1)/2
-j(n+1),
\endalign$$
 $m_0=\mult_{\{\O:\L\}}(1/p)$, $m_2=\mult_{\{\O:\L\}}(p)$,
and $\alpha_j(\L,\O)$ denotes the number of codimension $n-j$ totally
isotropic subspaces of $(\L\cap\O)/p(\L+\O)$.

Fix integral $\O\subseteq{1\over p}L$.  Decompose $\O$ as
$$\O={1\over p}\O_0\oplus\O_1\oplus p\O_2$$
where $\O_i\subseteq L$ with rank $r_i$, and $\O_0\oplus\O_1$ is
primitive in $L$ modulo $p$, meaning 
$(\O_0\oplus\O_1)\cap pL = p(\O_0\oplus\O_1)$, or equivalently,
$\dim(\overline{\O_0\oplus\O_1})$ in $L/pL$ is $r_0+r_1$.
(So $n=r_0+r_1+r_2$.)
Take $\L$ so that $p\O\subseteq \L \subseteq \left({1\over p}\O\cap L\right).$
Note that ${1\over p}\O\cap L=\O_0\oplus\O_1\oplus\O_2$.
Thus
$$\L=\O_0\oplus\left(\L_1\oplus p\L'_1\right)\oplus
\left(\L_2\oplus p\L'_2\oplus p^2\L''_2 \right)$$
where $\L_1\oplus\L_1'=\O_1$
and $\L_2\oplus\L'_2\oplus\L''_2=\O_2$.  Let $d_i=\rank \L_i$,
$d'_i=\rank\L'_i$, $d_2''=\rank\L_2''$.  Thus
$m_0=d_2$ and $m_2=r_0+d'_1+d''_2$.
Note that $\O_0$ is well-determined up to $p(\O_1\oplus\O_2)$,
and $\O_0\oplus\O_1$ is well-determined up to $p\O_2$.

Exactly as we did in [12], we can construct all these $\L$, simultaneously
constructing (and counting) all the subspaces of $(\L\cap\O)/p(\L+\O)$
counted by $\alpha_j(\L,\O)$.  In this way we find that the $\O$th coefficient
of $\t(L)|\widetilde T_j(p^2)$ is
$$\align
\widetilde c_j(\O)&=\sum_{\ell,d_2,t} (-1)^{j+r_0+d_2+t} p^E 
\varphi_{\ell}(\overline\O_1) \beta(r_2,x)
 p^{d_2(k-j+t)+d_2(d_2-1)/2}\beta(x,d_2)\\
&\qquad\cdot
\sum_{d_1'+d_2''=t} p^{d_2''(r_1-d_1'-\ell)}\beta(r_2-x,d_2'') 
\beta(r_1-\ell,d_1');
\endalign$$
here $x=j-r_0-\ell-t$ and
$$E=(k-n)(j-r_0-t)+(j-r_0-t)(j-r_0-t-1)/2+\ell(\ell+n-j-r_1+t).$$
Also, $\ell,t$ vary subject to $0\le \ell\le j-r_0$,
$0\le t\le j-r_0-\ell$.
By Lemma 5.1 (c) of [12], the sum on $d_1'+d_2''=t$ becomes
$\beta(r_1+r_2-x-\ell,t)=\beta(n-j+t,t)$.  Then by Lemma 4.2(a),
$$\sum_{d_2} (-1)^{d_2} p^{d_2(k-j+t)+d_2(d_2-1)/2}\beta(x,d_2)
= (-1)^{j-r_0-\ell-t} \mu(k-r_0-\ell-1,j-r_0-\ell-t).$$
Now, replacing $t$ by $j-r_0-\ell-t$,
and noting that 
$$\beta(m,r)\mu(m',r)= {\mu(m,r)\mu(m',r)\over\mu(r,r)}
=\beta(m',r)\mu(m,r),$$
 we get
$$\widetilde c_j(\O) = \sum_{\ell,t}(-1)^{\ell} p^{E'(\ell,t,\O)} 
\varphi_{\ell}(\overline\O_1) \beta(k-r_0-\ell-1,t)\mu(r_2,t)
\beta(n-r_0-\ell-t,j-r_0-\ell-t);$$
here $E'(\ell,t,\O)=\ell(k-r_0-r_1)+\ell(\ell-1)/2+t(k-n)+t(t-1)/2.$
This proves the proposition. $\square$
\enddemo

Next we extend Proposition 1.5 of [12].

\proclaim{Proposition 2.2} 
Suppose $1\le j\le k$ if $\chi(p)=1$, $1\le j<k$ if $\chi(p)=-1$.
Let $K_j$ vary over all lattices such that 
$pL\subseteq K_j\subseteq{1\over p}L$, 
$\mult_{\{L:K_j\}}\left({1\over p}\right) = \mult_{\{L:K_j\}}(p)=j,$
and $K_j\in\gen L$.
Then
$\sum_{K_j}\t(K_j;\tau)=\sum_{\O} b_j(\O)\e\{\O\tau\}$
where $\O$ varies over all even integral, (formally) rank $n$ sublattices
of ${1\over p}L$, and
$$\align
b_j(\O)
&= p^{(j-r_0)(j-r_0-1)/2} \sum_{\ell}p^{\ell(k-j-r_1+\ell)}
 \varphi_{\ell}(\overline\O_1)\\
&\qquad \cdot \delta(k-r_0-\ell-1,j-r_0-\ell)
\ \beta(k-r_0-r_1,j-r_0-\ell)
\endalign$$
if $\chi(p)=1$, 
$$\align
b_j(\O)&
= p^{(j-r_0)(j-r_0-1)/2} \sum_{\ell}(-1)^{\ell}
p^{\ell(k-j-r_1+\ell)} \varphi_{\ell}(\overline\O_1)\\
&\qquad \cdot \beta(k-r_0-\ell-1,j-r_0-\ell)
\ \delta(k-r_0-r_1,j-r_0-\ell)
\endalign$$
if $\chi(p)=-1$.
\endproclaim

\demo{Proof}
In Proposition 1.5 in [12] we showed that, for $p\not=2$,
$$b_j(\O)=p^{(j-r_0)(j-r_0-1)/2}\varphi_{j-r_0}(\overline\O_1^{\perp}\cap J)$$
where $L/pL=(\overline\O_0\oplus\overline\O_0')\perp J$,
$\overline\O_0\oplus\overline\O_0\simeq\H^{r_0}$, $\overline\O_1\subseteq J$.
Using the results of \S5, this argument is valid for $p=2$.
Since $L/pL$ is regular, so is $J$, and $J$ is hyperbolic if and only if
$L/pL$ is.  Also, we necessarily have $r_0\le k$ if $\chi(p)=1$
(and hence $L/pL$ is hyperbolic), and $r_0<k$ otherwise.
Decompose $\overline\O_1=\O_1/p\O_1$ as $R\perp W$ where $R=\rad\overline\O_1$;
so $W$ is regular.  Then $J=(R\oplus R')\perp W\perp W'$ where
$R\oplus R'\simeq\H^r$, $r=\dim R$, and $W'$ regular 
(we use 42:4 of [8] when $p\not=0$, and Proposition 5.2 when $p=2$).

(a) Say $\chi(p)=1$.  So $L/pL$ and $J$ are hyperbolic.  Hence if
$W$ is hyperbolic, $W'$ must be as well.  If $W\simeq\H^d\perp\A$, $\A$ an
anisotropic plane, then we must have $W'\simeq\H^{d'}\perp\A$ (some $d,d'$).
If $W$ has odd dimension, then so does $W'$.

Here $\dim J=2(k-r_0)$, $\dim W=r_1-r$, $\dim W'=2(k-r_0)-r_1-r$, and
$(R\perp W)^{\perp}=(R\perp W')$ (which has dimension $2(k-r_0)-r_1$).
So $W'=W''\perp\H^{k-r_0-r_1}$ and $R\perp W'=R\perp W''\perp\H^{k-r_0-r_1}$
where $W''$ is regular of dimension $r_1-r$, with $W''$ hyperbolic if
and only if $W$ is.  (See the discussion at the beginning of \S4 to
make sense of $W\perp\H^t$ when $t<0$.)  Consequently, recalling the 
formulas for $\varphi_{\ell}(*)$ from the beginning of this section,
$$
\varphi_{\ell}(\overline\O_1\cap J)
=\varphi_{\ell}(R\perp W''\perp\H^{k-r_0-r_1})
=\varphi_{\ell}(\overline\O_1\perp\H^{k-r_0-r_1}).
$$
We now apply our Reduction Lemma (Lemma 4.1) to obtain the result.

(b) Now say $\chi(p)=-1$.  Thus $L/pL\simeq\H^{k-1}\perp\A$, and 
$J\simeq\H^{k-r_0-1}\perp\A$.  Thus with analysis virtually identical to
that used above, we find
$$\varphi_{\ell}(\overline\O_1^{\perp}\cap J)
= \varphi_{\ell}(\overline\O_1\perp\H^{k-r_0-r_1-1}\perp\A).$$
Now apply the Reduction Lemma (Lemma 4.1).  
$\square$
\enddemo

These last two results allow us to prove the following.

\proclaim{Theorem 2.3}  Say $1\le j\le k$ if $\chi(p)=1$, 
$1\le j<k$ if $\chi(p)=-1$.
Let $K_j$ be as in Proposition 2.2.  Set
$$u_q(j)=(-1)^q p^{q(q-1)/2}\beta(n-j+q,q),\qquad
T'_j(p^2) = \sum_{0\le q\le j} u_q(j) \widetilde T_{j-q}(p^2),$$
$$v_q(j)= 
\cases (-1)^q \beta(k-n+q-1,q)\delta(k-j+q-1,q)&\text{if $\chi(p)=1$,}\\
(-1)^q \delta(k-n+q-1,q)\beta(k-j+q-1,q)&\text{if $\chi(p)=-1$.}
\endcases$$
 Then
$$\t(L)|T'_j(p^2) = 
\sum_{0\le q\le j} v_q(j)\left(\sum_{K_{j-q}}\t(K_{j-q})\right)$$
where $K_{j-q}$ varies subject to $pL\subseteq K_{j-q}\subseteq{1\over p}L$,
$$\mult_{\{L:K_{j-q}\}}(1/p) = \mult_{\{L:K_{j-q}\}}(p) = j-q,$$ and
$K_{j-q}\in\gen L$.
\endproclaim

\demo{Proof}
When $\chi(p)=1$, this is proved in Theorem 1.2 of [12].  So
suppose $\chi(p)=-1$.
We show that
$$\align
&\sum_q u_q(j) \widetilde c_{j-q}(\O)\\
&\qquad= \sum_{\ell} (-1)^{\ell} p^{E''_j(\ell,\O)}
\varphi_{\ell}(\overline\O_1)
\beta(k-r_0-\ell-1,j-\ell-r_0) \mu(r_2,j-\ell-r_0) \\
&\qquad= \sum_q v_q(j) b_{j-q}(\O),
\endalign$$
where
$$\align
E''_j(\ell,\O)&=E'(\ell,j-r_0-\ell,\O) \\
&=(j-r_0)(j-r_0-1)/2+(k-n)(j-r_0)-\ell(j-r_0-\ell-r_2),
\endalign$$
$E'(\ell,t,\O)$ as defined in Proposition 2.1.

Using the formula in Proposition 2.2 for $\widetilde
c_{j-q}(\O)$, we have
$$\align
\sum_q u_q(j) \widetilde c_{j-q}(\O)
& = \sum_{\ell,t} (-1)^{\ell} p^{E'(\ell,t,\O)} \varphi_{\ell}(\overline\O_1)
\beta(k-r_0-\ell-1,t)\mu(r_2,t)\\
&\qquad\cdot \sum_q u_q(j) \beta(n-r_0-\ell-t,j-r_0-\ell-t).
\endalign$$
We claim the sum on $q$ is 0 unless $t=j-r_0-\ell.$  Using Lemma 5.1 (b)
of [12] with $r=n-r_0-\ell-t$, $m=j-r_0-\ell-t$, $m'=q$, we have
$$\align
&\sum_q u_q(j) \beta(n-r_0-\ell-t,j-q-r_0-\ell-t) \\
&\qquad = \beta(n-r_0-\ell-t,j-r_0-\ell-t)
\sum_q (-1)^q p^{q(q-1)/2} \beta(j-r_0-\ell-t,q).
\endalign$$
By Lemma 4.2 (a), the latter sum on $q$ is 0 provided $t<j-r_0-\ell$;
when $t=j-r_0-\ell$, the sum on $q$ is 1.
Thus 
$$\sum_q u_q(j) \widetilde c_{j-q}(\O) 
=\sum_{\ell} (-1)^{\ell} p^{E''_j(\ell,\O)} \varphi_{\ell}(\overline\O_1)
\beta(k-r_0-\ell-1,j-r_0-\ell) \mu(r_2,j-r_0-\ell).$$

On the other hand, using Lemma 5.1 (b) of [12], and with
$S(m)$ defined as in Lemma 4.2 (b),
we have
$$\align
\sum_q v_q(j) b_{j-q}(\O)
&= \sum_{\ell} (-1)^{\ell} p^{(j-r_0)(j-r_0-1)/2 + \ell(k-r_1-j+\ell)}
\varphi_{\ell}(\overline\O_1)\\
&\qquad \cdot \beta(k-r_0-\ell-1,j-r_0-\ell) S(j-r_0-\ell).
\endalign$$
Applying Lemma 4.2 (b) completes the proof. $\square$
\enddemo

We say $K$ lies in the genus of $L$, denoted
$K\in\gen L$, if for all primes $q$, $\Z_q K\simeq\Z_q L$.  With
$o(K)$ the order of the orthogonal group of $K$, we set
$$\t(\gen L)=\sum_{\cls K}{1\over o(K)}\t(K)$$
where $\cls K$ runs over all isometry classes in the genus of $L$.
(Note: Sometimes people use $\t(\gen L)$ to refer to the normalized
average 
$${1\over\mass L}\sum_{\cls K} {1\over o(K)} \t(K)$$ where
$\mass L=\sum_{\cls K} {1\over o(K)}.$)

\proclaim{Corollary 2.4}  Suppose
$j\le k$ when $\chi(p)=1$, and $j<k$ when $\chi(p)=-1$,
$$\t(\gen L)|T'_j(p^2)=\lambda_j(p^2) \t(\gen L)$$ where
$$\align&
\lambda_j(p^2)=\\
&\qquad\cases p^{j(k-n)+j(j-1)/2}\beta(n,j)(p^{k-1}+1)(p^{k-2}+1)
\cdots(p^{k-j}+1)
&\text{if $\chi(p)=1$,}\\
p^{j(k-n)+j(j-1)/2}\beta(n,j)\mu(k-1,j)
&\text{if $\chi(p)=-1$.}\endcases
\endalign$$
\endproclaim

\demo{Proof}  
The case $\chi(p)=1$ was treated in Corollary 1.3 of [12].  So suppose
$\chi(p)=-1$.
We average across the identity of the theorem,
getting
$$\t(\gen L)|T'_j(p^2)
= \sum_q v_q(j)\left( \sum_{\cls L'}{1\over o(L')}
\sum_{K_{j-q}} \t(K_{j-q}) \right).$$
(Here $\cls L'$ varies over $\gen L$.)  As we argued in the proof
of Corollary 1.3 of [12],
we have
$$\sum_{\cls L'}{1\over o(L')} \sum_{K_m'}\t(K_m')
= \sum_{\cls K'} p^{m(m-1)/2}\varphi_m(L/pL)\cdot {1\over o(K')}\t(K').$$
Here the lattices $K_m'$, $pL'\subseteq K_m'\subseteq{1\over p}L'$,
 vary as in Proposition 2.2, 
and $\cls K'$ varies over $\gen L$.
Thus $\sum_{\cls K'}{1\over o(K')}\t(K')=\t(\gen L)$, and so
$$\t(\gen L)|T'_j(p^2)
=\lambda_j(p^2) \t(\gen L)$$
where 
$\lambda_j(p^2)=\sum_q v_q(j) p^{(j-q)(j-q-1)/2} \varphi_{j-q}(L/pL).$
Since $\chi(p)=-1$, we know $L/pL\simeq\H^{k-1}\perp\A$; thus 
using the formula for $\varphi_{\ell}(U)$ presented at
the beginning of this section, we have
$$\varphi_m(L/pL)=\prod_{i=0}^{m-1} {(p^{k-i}+1)(p^{k-i-1}-1)\over(p^{m-i}-1)}
= \delta(k,m) \beta(k-1,m).$$
Thus again using Lemma 5.1 (b) of [12], we get 
$\lambda_j(p^2)= p^{j(j-1)/2}\beta(k-1,j)S(j)$ where
$$S(j)=
\sum_{q=0}^j p^{q(q+1)/2-qj} \delta(k-n+q-1,q) \delta(k,j-q) \beta(j,q).$$
To evaluate $S(j)$, we use the identity
$$\beta(j,q)=\beta(j-1,q)+p^{j-q}\beta(j-1,q-1)$$
to split the sum defining $S(j)$ into a sum on $0\le q<j$ and on
$1\le q\le j$.  Then we replace $q$ by $q+1$ in the second sum.
Arguing by induction on $d$ with the hypothesis
$S(j)=p^{d(k-n)}\mu(n-j+d,d) S(j-d)$ now easily gives us the value of
$\lambda_j(p^2)$, as claimed.  $\square$
\enddemo

\head{\S3. Vanishing of theta series under Hecke operators,
 and  $\t(\gen L)|T(p)^2$}
\endhead
As in the preceeding section, 
$L$ is an even integral, rank $2k$ lattice with positive definite quadratic
form of level $N$; we
fix a prime $p$, $p\nmid N$.

Although Proposition 2.1 is valid for all values
of $j\le n\le 2k$,
the geometry of $L/pL$ presents an obstruction to extending Proposition 2.2
for $j>k$ when $\chi(p)=1$, and for $j\ge k$ when $\chi(p)=-1$.
 However, given even integral 
$$\O={1\over p}\O_0\oplus\O_1\oplus p\O_2\subseteq {1\over p}L$$
with $\O_0\oplus\O_1$ primitive in $L$ modulo $p$
(meaning $(\O_0\oplus\O_1)\cap pL=p(\O_0\oplus\O_1)$),
we necessarily have $r_0=\rank\O_0\le k$ if $\chi(p)=1$, $r_0<k$ if
$\chi(p)=-1$.
Thus all the $\O\subseteq{1\over p}L$ that arise when describing
$\t(L)|T'_j(p^2)$ for ``large'' $j$ (i.e. $j>k$ when $\chi(p)=1$,
$j\ge k$ when $\chi(p)=-1$)
have already been considered when describing
$\t(L)|T'_j(p^2)$ for small $j$.  In fact, we find the following.

\proclaim{Proposition 3.1}  
For $q\ge 0$, $a\ge 1$, set
$$w_q(a)=
\cases(-1)^q p^{q(q+1)/2} \beta(a+q-1,q) \beta(n-k+q,a+q)
&\text{if $\chi(p)=1$,}\\
(-1)^q p^{q(q+1)/2} \beta(a+q-1,q) \beta(n-k+1+q,a+q)
&\text{if $\chi(p)=-1$.}
\endcases$$
Then
$$\t(L)|\widetilde T_{k+a}(p^2) = 
\t(L)|\left( \sum_{0\le q\le k} w_q(a)\widetilde T_{k-q}(p^2)\right)$$
when $\chi(p)=1$, and
$$\t(L)|\widetilde T_{k-1+a}(p^2) = 
\t(L)|\left( \sum_{0\le q\le k} w_q(a)\widetilde T_{k-1-q}(p^2)\right)$$
when $\chi(p)=-1$.
\endproclaim

\demo{Proof}  The proofs for the cases $\chi(p)=1$ and $\chi(p)=-1$ are
virtually identical, so we present only the case when $\chi(p)=1$.

We prove 
$$\widetilde c_{k+a}(\O) = 
\sum_{0\le q\le k} w_q(a)\widetilde c_{k-q}(\O).$$
In our formula for $\widetilde c_j(\O)$ given in Proposition 2.1, 
only one term
is dependent on $j$.  Thus proving our claim reduces to proving that
for $0\le \ell\le k-r_0$, $0\le t\le k-r_0-\ell$, and $x=k-r_0-\ell-t$,
$$\sum_{0\le q\le x} w_q(a) \beta(n-k+x,x-q)=\beta(n-k+x,x+a).$$
By Lemma 5.1(b) of [12],
$$\beta(n-k+x,x-q)\beta(n-k+q,a+q) = \beta(n-k+x,x+a)\beta(x+a,x-q).$$
This identity together with Lemma 4.2(c) establishes the claim.
$\square$
\enddemo

We want to state our main results in terms of $T'_j(p^2)$ rather than
$\widetilde T_j(p^2)$.  To aid with this, we have:

\proclaim{Proposition 3.2}  Let $T'_j(p^2)$ be defined as in Theorem 2.3.
For $r\ge 1$,
$\widetilde T_r(p^2) = \sum_{0\le q\le r} \beta(n-q,r-q) T'_q.$
\endproclaim

\demo{Proof}  We evaluate the right-hand side expression by first
replacing $q$ by $r-q$, then substituting for $T'_{r-q}$ in terms of
$\widetilde T_{r-q-i}$, $0\le i\le r-q$.  Then we replace $i$ by $r-q-i$.
This gives us a sum over $0\le q\le r$, $0\le i\le r-q$, or equivalently,
changing the order of summation, a sum over $0\le i\le r$,
$0\le q\le r-i$.  Finally, we evaluate the sum on $q$ using Lemma 4.2 (c).
$\square$
\enddemo

Now we can prove:

\proclaim{Theorem 3.3}  
Say $1\le a\le k$ if $\chi(p)=1$, $1\le a\le k+1$ if $\chi(p)=-1$.
We have $\t(L)|T_{k+a}'(p^2) = 0$
when $\chi(p)=1$, and
$\t(L)|T_{k-1+a}'(p^2) = 0$ when $\chi(p)=-1$.
\endproclaim

\demo{Proof} Say $\chi(p)=1$.  By Proposition 3.1, we see
$$\t(L)|\widetilde T_{k+a} 
= \sum_{0\le q\le k} w_q(a) \t(L)|\widetilde T_{k-q}.$$ 
Using also Proposition 3.2,
we find
$$
\t(L)|T_{k+a}'
= \t(L) | \sum_{q,\ell} w_q(a) \beta(n-\ell,k-q-\ell) T_{\ell}' 
 - \sum_r \beta(n-r,k+a-r) T_r'.
$$
Here $0\le r<k+a$; also, $0\le q\le k$, $0\le\ell\le k-q$, or equivalently,
$0\le \ell\le k$, $0\le q\le k-\ell$.
Using first Lemma 5.1(b) of [12] and then Lemma 4.2(c), we find
$$\align
& \sum_{0\le q\le k-\ell} w_q(a) \beta(n-\ell,k-\ell-q) \\
&\qquad = \beta(n-\ell,k-\ell+a) \sum_q (-1)^q p^{q(q+1)/2}
\beta(a+q-1,q) \beta(k-\ell+a,k-\ell-q) \\
&\qquad = \beta(n-\ell,k-\ell+a).
\endalign$$
Thus $\t(L)|T'_{k+1} = 0$, and 
$\t(L)|T'_{k+a} = -\t(L)|\sum_r \beta(n-r,k+a-r) T'_r$
where $k<r<k+a$.  Hence by induction on $a$,
$\t(L)|T'_{k+a} = 0$ for all $a\ge 1$.

When $\chi(p)=-1$ the proof is virtually identical. $\square$
\enddemo

In Proposition 1.4 of [12] we showed that if $n\le k$ then 
$$\t(\gen L)|T(p)=\epsilon(k-j,n)\t(\gen L')
=\delta(k-1,n) \t(\gen L')$$ when $\chi(p)=1,$
and $\gen L=\gen L'$ when $\chi=1$.  
(In fact, for primes $q\not=p$, $\Z_q L'\simeq\Z_q L$, and
$\Z_p L'\simeq\Z_p L$.  If $\left(p\over q\right)=1$ for all primes
$q|N$, then $\gen L'=\gen L$ and so $\t(\gen L)$ is an eigenform 
for $T(p)$.)
Using Lemma 4.1, 
we can show this holds for all
$n\le 2k$.  We can also extend this result to include $T(p)^2$, $\chi(p)=\pm1$.

\proclaim{Theorem 3.4}  Say $n\le 2k$.
\item{(a)} 
If $\chi(p)=1$ then $\t(\gen L)|T(p) = \delta(k-1,n) \t(\gen K^{1/p})$
where $pL\subseteq K\subseteq L$ with
$\Z_p K^{1/p} \simeq \Z_p L$.  When $\chi=1$, $\gen K=\gen L$ and so
$\t(\gen L)$ is an eigenform for $T(p)$.
\item{(b)} If $\chi(p)=1$ then $\t(\gen L)|T(p)^2 = 
\big( \delta(k-1,n) \big)^2 \t(\gen L).$
\item{(c)} If $\chi(p)=-1$ then $\t(\gen L)|T(p)^2 = 
\big( \mu(k-1,n) \big)^2 \t(\gen L).$
\endproclaim

\demo{Proof}  (a) Lemma 4.1 extends Lemma 1.6 of [12],
which allows us to prove this
result for $n\le 2k$, just as we did for $n\le k$ in [12].

The proofs of (b) and (c) are virtually identical, so we prove (c).
By a straightforward extension of Proposition 5.1 of [5]
with $\chi(p)=-1$, we have
$$T(p)^2 = \sum _{0\le j\le n} (-1)^{n-j} p^{k(n-j)+j(j+1)/2 - n(n+1)/2}
\widetilde T_j(p^2).$$
Replacing $j$ by $n-j$ and then using Proposition 3.2,
$$T(p)^2 = \sum_{j,r} (-1)^j p^{j(j-1)/2+j(k-n)} \beta(n-r,j) T'_r,$$
where $0\le j\le n$, $0\le r\le n-j$, or equivalently, $0\le r\le n$,
$0\le j\le n-r$.  Using that $\beta(n-r,j)=p^j\beta(n-r-1,j)+\beta(n-r-1,j-1)$,
we split the sum on $j$ into a sum on $0\le j<n-r$ and a sum on $0<j\le n-r$.
Then we replace $j$ by $j+1$ in the second sum and simplify to get
$$T(p)^2 = \sum_{0\le r\le n} (-1)^{n-r} \mu(k-1-r,n-r) T'_r(p^2).$$
Thus $\t(\gen L)|T(p)^2 = \lambda \t(\gen L)$ where
$$\align
\lambda 
&= (-1)^n \sum_{0\le r\le n} (-1)^r p^{r(r-1)/2 +r(k-n)}
\mu(k-1-r,n-r) \mu(k-1,r) \beta(n,r)\\
&= (-1)^n \mu(k-1,n) S(n,k-n)
\endalign$$
where $S(n,k-n)$ is defined and evaluated in Lemma 4.2(a).
We quickly find $\lambda = \big( \mu(k-1,n)\big)^2.$

A similar argument shows that when $\chi(p)=1$,
$$\t(\gen L)|T(p)^2 = \big(\delta(k-1,n)\big)^2 \t(\gen L).\qquad
\square$$
\enddemo

\head{\S4. Lemmas on quadratic spaces over $\Z/p\Z$}\endhead
In this section we rely on \S42 and \S62 of [8] for results on quadratic 
spaces over a field $\F$ with $p=\char \F\not=2$.  
One can deduce from
\S93 of [8] the corresponding results when $p=2$; for completeness
we present the results we need in \S 5.

Our first goal of this section is to prove our Reduction Lemma
(Lemma 4.1), critical
to the proof of our main
 theorem.  This lemma focuses on a general formula that 
counts totally isotropic subspaces of fixed degree within any given 
quadratic space over $\Z/p\Z$.
In [12] we proved a formula to allow us to ``cancel'' hyperbolic
planes from the expression $\varphi_{\ell}(U\perp\H^t)$ provided $t\ge 0$.
In Lemma 4.1 we extend this result to allow $t<0$, and to allow us  also to cancel
anisotropic planes.

Note that if $U\simeq W\perp\H^t$ for some $t\ge 0$, then $U\perp\H^{-t}$ is
meaningful.  
Also, by 42:16 [8] and Theorem 5.7,
if $U\simeq W\perp\H^t\simeq W'\perp\H^t$ then $W\simeq W'$
(and so $U\perp\H^{-t}$ is well-defined up to isometry whenever the
expression is meaningful).
  Similarly, if $U\simeq W\perp\H^t\perp\A$, $t\ge 0$, then
$U\perp\H^{-t}\perp\A^{-1}$ is meaningful, and since 
$\H\perp\H\simeq\A\perp\A$,
we can write this as $U\perp\H^{-t-2}\perp\A$.  
Also, by 
if $U\simeq W\perp\H^t\perp\A\simeq W'\perp\H^t\perp\A$ then $W\simeq W'$
(see 42:16 of [8] when $p$ is odd, and Theorems 5.7 and 5.9 when $p=2$).
Thus $U\perp\H^{-t}\perp\A^{-1}$ is well-defined up to isometry whenever
the expression is meaningful.

\proclaim{Lemma 4.1 (Reduction Lemma)} Let $U$ be a dimension $d$ space
over $\Z/p\Z$, $\ell\ge0$ and $t\in\Z$ so that $U\perp\H^t$
(resp. $U\perp\H^t\perp\A$) is defined.
\item{(a)}
$\varphi_{\ell}(U\perp\H^t) = \sum_{r=0}^{\ell} p^{r(t-\ell+r)}
\delta(d-1+t-r,\ell-r) \beta(t,\ell-r) \varphi_r(U).$
\item{(b)} 
$\varphi_{\ell}(U\perp\H^t\perp \A)
= \sum_{r=0}^{\ell} (-1)^r p^{r(t+1-\ell+r)}
\beta(d+t-r,\ell-r) \delta(t+1,\ell-r) \varphi_r(U).$
\endproclaim

\demo{Proof}  
(a)
  We established this for $t\ge 0$ in Lemma 5.1 [12].  Now fix $t>0$,
and say $W\simeq U\perp\H^{-t}$ (in other words, $W$ is a lattice so that
$U\simeq W\perp\H^t$).  Then using that $\delta(m,r')\delta(m-r',r)=
\delta(m,r+r')$, we get
$$\align
& \sum_{r=0}^{\ell} p^{r(-t-\ell+r)}
\delta(d-1-t-r,\ell-r) \beta(-t,\ell-r) \varphi_r(U) \\
&\qquad = \sum_{r=0}^{\ell} p^{(\ell-r)(-t-r)}
\delta(d-1+t-\ell+r,r) \beta(-t,r) \varphi_{\ell-r}(W\perp\H^t) \\
&\qquad = \sum_{r,q} p^{-t\ell+q(t-\ell+q)}\delta(d-1-t-q,\ell-q)
\varphi_q(W) S(\ell-q);
\endalign$$
here $0\le r\le \ell$, $0\le q\le \ell-r$, or equivalently,
$0\le q\le \ell$, $0\le r\le \ell-q$, and
$$S(m)=\sum_{r=0}^m p^{r(r+t-m)}\beta(-t,r)\beta(t,m-r).$$
By definition, we see $\beta(-t,r)=(-1)^r p^{-rt-r(r-1)/2} \beta(t+r-1,r).$
Using this identity and Lemma 4.2 (c) (with $y=0, a=t$), 
we see $S(m)=1$ if $m=0$, and 0 otherwise.
Thus our final sum on $q,r$ simplifies to be $\varphi_{\ell}(W)
=\varphi_{\ell}(U\perp\H^{-t}),$ proving part (a) of the lemma.

(b)  We first establish the lemma for $t=0$.  
If $d=0$ the claim is trivial.  If $d>0$ and $U$ is anisotropic
then either $U\simeq <\varepsilon>$ ($\varepsilon\not=0$) or $U\simeq\A$.
The claim is easily established in either case, noting that
$<\varepsilon>\perp\A\simeq\H\perp<-\delta\varepsilon>$ where
$\left({\delta\over p}\right)=-1$, and $\A\perp\A\simeq\H\perp\H$.
Also recall that $\beta(d-k,\ell-k)=0$ if $\ell>d\ge 0$.

Now suppose $U$ is isotropic.  Then either $U\simeq W\perp<0>$ or
$U\simeq W\perp\H$.  When $d=1$ or 2, i.e. $U\simeq<0>, <0,0>$ or $\H$,
the claim is easily established.  

We argue by induction on $d$,
so now we suppose $d>2$ and 
that the lemma holds for spaces of dimension less than $d$.

First suppose $U\simeq W\perp<0>.$  Then a dimension $r$ totally
isotropic subspace of $U$ projects onto a dimension $r$ or $r-1$
totally isotropic subspace of $W$.  Thus
$$\varphi_r(U)=p^r\varphi_r(W) + \varphi_{r-1}(W),$$
and similarly,
$$\varphi_{\ell}(U\perp\A)
=p^{\ell} \varphi_{\ell}(W\perp\A) + \varphi_{\ell-1}(W\perp\A).$$
In any case, $\dim W=d-1$ so the lemma holds for $\varphi_*(W\perp\A)$.
Using this, we get
$$\align
\varphi_{\ell}(U\perp\A)
&= (-1)^{\ell} p^{2\ell} \varphi_{\ell}(W) \\
&\qquad + \sum_{r=0}^{\ell-1}(-1)^r p^{r(r-\ell+2)}\varphi_r(W) 
\delta(1,\ell-r-1)\\
&\qquad\qquad \cdot \big[\beta(d-r,\ell-r)+p^2 \beta(d-r-1,\ell-r) \big]. \\
\endalign$$
On the other hand,
$$\align
& \sum_{r=0}^{\ell} (-1)^r p^{r(r-\ell+1)} \delta(1,\ell-r) 
\beta(d-r,\ell-r) \varphi_r(U) \\
&\qquad = (-1)^{\ell} p^{2\ell} \varphi_{\ell}(W) \\
&\qquad \qquad + \sum_{r=0}^{\ell-1}(-1)^r p^{r(r-\ell+2)}\varphi_r(W) \\
&\qquad\qquad \cdot \delta(1,\ell-r-1)
\big[\beta(d-r,\ell-r)+p^2 \beta(d-r-1,\ell-r) \big]. \\
\endalign$$
Thus the lemma holds for $U$ of dimension $d$, $U\simeq W\perp<0>.$

Now say $U\simeq W\perp\H$.  Then by Lemma 4.1 [12] we have
$$\varphi_r(W\perp\H) = p^r \varphi_r(W) + (p^{d-1-r}+1)\varphi_{r-1}(W)$$
and so
$\varphi_{\ell}(U\perp\A)
=p^{\ell} \varphi_{\ell}(W\perp\A) 
 + (p^{d+1-r}+1) \varphi_{\ell-1}(W\perp\A).$
Applying our induction hypothesis to $\varphi_*(W\perp\A)$ and simplifying,
$$\align
\varphi_{\ell}(U\perp\A)
&= \sum_{r=0}^{\ell}(-1)^r p^{r(r-\ell+2)}\varphi_r(W) 
\delta(2,\ell-r) \beta(d-r-1,\ell-r)
\endalign$$
On the other hand, applying Lemma 1.6 [12] on $\varphi_r(U)=
\varphi_r(W\perp\H)$, we find
$$\align
& \sum_{r=0}^{\ell} (-1)^r p^{r(r-\ell+1)} \delta(1,\ell-r) 
\beta(d-r,\ell-r) \varphi_r(U) \\
&\qquad = 
\sum_{r=0}^{\ell}(-1)^r p^{r(r-\ell+2)}\varphi_r(W) 
\delta(2,\ell-r) \beta(d-r-1,\ell-r)
\endalign$$
This proves the lemma for $U\simeq W\perp\H$, $d=\dim U$.  Induction on $d$ now
shows the lemma holds for $t=0$ and all $\ell,d\ge 0$.

Now fix $t\ge 0$, and say the lemma holds for this $t$ and all $\ell,d\ge 0$.
By Lemma 4.1 [12] and then the induction hypothesis, we get
$$\align
& \varphi_{\ell}(U\perp\H^{t+1}\perp\A) \\
&= p^{\ell}\varphi_{\ell}(U\perp\H^t\perp\A)
+(p^{d+2t+3-\ell}+1) \varphi_{\ell-1}(U\perp\H^t\perp\A)\\
&=p^{\ell}\sum_{r=0}^{\ell} (-1)^r p^{r(t+1-\ell+r)} \varphi_r(U)
\delta(t+1,\ell-r)
\beta(d+t-r,\ell-r) \\
&\qquad + (p^{d+2t+3-\ell}+1) \sum_{r=1}^{\ell-1} (-1)^r p^{r(t+2-\ell+r)}
\varphi_r(U) \\
&\qquad\qquad \cdot \delta(t+1,\ell-r-1)\beta(d+t-r,\ell-r-1) \\
& = \sum_{r=0}^{\ell} (-1)^r p^{r(t+2-\ell+r)}\delta(t+2,\ell-r)
\beta(d+t+1-r,\ell-r)\varphi_r(U),
\endalign$$
as claimed.  Induction on $t$ now shows the lemma holds for all
$t,\ell,d\ge 0$.

To show the formula holds for $\varphi_{\ell}(U\perp\H^{-t}\perp\A)$ when
$t>1$, we set $W=U\perp\H^{-t}\perp\A$.  So
$W\simeq U\perp\H^{t-2}\perp\A$ (recall $\H\perp\H\simeq \A\perp\A,$
so $\H^t\perp\A^{-1}\simeq \H^{t-2}\perp\A$).
Then we proceed 
just as we did to verify our formula for $\varphi(U\perp\H^{-t})$,
noting that $\delta(-t+1,s)=p^{(-t+1)s-s(s-1)/2} \delta(s+t-2,s)$
and using Lemma 4.2(b).

Finally, to prove our formula for $\varphi_{\ell}(U\perp\H^{-1}\perp\A)$,
we proceed by induction on $\ell\ge 0$.  The formula clearly holds for 
$\ell=0$.  So assume the formula holds for 
$\varphi_{\ell-1}(U\perp\H^{-1}\perp\A)$.  By what we have already proved,
$$\varphi_{\ell}((U\perp\H^{-1}\perp\A)\perp\H)
= p^{\ell} \varphi_{\ell}(U\perp\H^{-1}\perp\A) + 
(p^{d+1-\ell}+1) \varphi_{\ell-1}(U\perp\H^{-1}\perp\A).$$
 Also, we know
$\varphi_{\ell}((U\perp\H^{-1}\perp\A)\perp\H) =
\varphi_{\ell}(U\perp\A).$
Thus using our formula with $t=0$, our induction hypothesis, and
straightforward simplification,
we get
$$\align
&p^{-\ell}\varphi_{\ell}(U\perp\H^{-1}\perp\A)\\
&\qquad
= -(p^{d+\ell}+1)\sum_{r=0}^{\ell-1} (-1)^r p^{r(1-\ell+r)}
\delta(0,\ell-r-1) \beta(d-r-1,\ell-r-1) \varphi_r(U) \\
&\qquad + \sum_{\ell=0}^{\ell} (-1)^r p^{r(1-\ell+r)}
\delta(1,\ell-r) \beta(d-r,\ell-r) \varphi_r(U) \\
&\qquad
= \sum_{r=0}^{\ell-1} (-1)^r p^{r(r-\ell)+\ell} \delta(0,\ell-r)
\beta(d-1-r,\ell-r) \varphi_r(U).
\endalign$$
This finishes proving the lemma. $\square$
\enddemo

In the following lemma we collect some identities we found quite useful.

\proclaim{Lemma 4.2} Fix $a,m,y\in\Z$ with $a,m>0$.
\item{(a)} Set $S(m,y)=\sum_{0\le q\le m} (-1)^q p^{q(q-1)/2+qy}\beta(m,q).$
Then $S(m,y)=(-1)^m \mu(y+m-1,m).$
\item{(b)}  Let $S(m)=\sum_{0\le q\le m} (-1)^q p^{q(q+1)/2+q(y-m)}
\delta(a-1+q,q) \delta(a+y,m-q) \beta(m,q).$
Then $S(m)=(-1)^m \mu(y,m).$
\item{(c)}  Let $S(m)=\sum_{0\le q\le m}(-1)^q p^{q(q+1)/2+q(y-m)}
\mu(a-1+q,q)\mu(a+y,m-q)\beta(m,q).$
Then
$$\sum_{0\le q\le m} (-1)^q p^{q(q+1)/2+q(y-m)}
\beta(a-1+q,q) \beta(a+y,m-q)={1\over \mu(m,m)} S(m),$$
and $S(m)=\mu(y,m).$
\item{(d)} Let $S(m)=\sum_{0\le q\le m} (-1)^q p^{q(q+1)/2}
\mu(a+m-q,m-q) \mu(b,q) \beta(m,q).$
Then
$$\align
S(m) 
&=\sum_{0\le q\le m} (-1)^q p^{q(q+1)/2}
\delta(a+m-q,m-q) \delta(b,q) \beta(m,q) \\
&= (-1)^m p^{am+m(m+1)/2} \mu(b-a-1,m).
\endalign$$
\endproclaim

\demo{Proof}  The method of proof is virtually identical for all these
claims, so we prove here only one of the claims made in (d) and we comment on
how to adapt this proof to prove the other claims in the lemma.

Using that $\beta(m,q)=p^q \beta(m-1,q)+\beta(m-1,q-1)$ when $0<q<m$,
we can separate the sum defining $S(m)$ into a sum on $0\le q<m$ and
a sum on $0<q\le m$.  Then replacing $q$ by $q+1$ in the second sum,
we find
$$\align
S(m)&= \sum_{0\le q<m} (-1)^q p^{q(q-1)/2} \mu(a+m-q-1,m-q-1) \mu(b,q)
\beta(m-1,q)\\
&\qquad \cdot \left[ p^q(p^{a+m-q}-1) - p^q(p^{b-q}-1)\right] \\
&= - p^{a+m} (p^{b-a-m}-1) S(m-1).
\endalign$$
Arguing by induction on $d$ with the hypothesis
$$S(m)=(-1)^d p^{d(a+m)-d(d-1)/2} \mu(b-a-m+d-1,d) S(m-d),$$
we show $S(m)=(-1)^m p^{am+m(m+1)/2}\mu(b-a-1,m)$, as claimed.

A virtually identical argument (with $\mu$ replaced by
$\delta$) proves the other claim in (d).

To prove (a), (b), and 
(c),  we split the sum as above, and, as above, we then shift
the variable in one sum and simplify.
For (a) we use the induction hypothesis
$$S(m,y)=(-1)^d \mu(y+d-1,d) S(m-d,y+d);$$
for (b) we use $S(m)=(-1)^d \mu(d+y-m,d) S(m-d).$
To prove (c) we use the identity
$$\align
&\beta(a-1+q,q) \beta(a+y,m-q) \\
&\qquad = {1\over \mu(m,m)}\cdot
\mu(a-1+q,q) \mu(a+y,m-q) {\mu(m,m)\over \mu(q,q)\mu(m-q,m-q)} \\
&\qquad = {1\over \mu(m,m)}\cdot
\mu(a-1+q,q) \mu(a+y,m-q) \beta(m,q), \qquad 
\endalign$$
and the induction hypothesis that
$S(m)=\big(\mu(y-m+d,d)/\mu(m,d)\big) S(m-d).$
$\square$
\enddemo

\head{\S5. Quadratic spaces over fields of characteristic 2}
\endhead

Let $\F$ be a characteristic 2 field of order $q$, $V$ a finite-dimensional
space over $\F$ with quadratic form $Q'$ and symmetric bilinear form $B$
such that $Q'(x+y)=Q'(x)+Q'(y)+B(x,y)$.  (This is the relative scaling between
$Q'$ and $B$ used in [2].)  Note that $B(x,x)=2Q'(x)$, so every vector
is orthogonal to itself.

Define the radical of $V$ to be
$$\rad V=\{x\in V:\ Q'(x)=0 \text{ and }B(x,V)=0\ \}.$$
We say $V$ is regular if $\rad V=\{0\}$.  
Clearly $\rad V$ is a subspace of $V$, and $V=\rad V\perp U$ for some
regular subspace $U$.  While $U$ is not uniquely determined, if
$V=\rad\perp U=\rad V\perp U'$ then there exist bases
$(u_1,\ldots,u_d)$, $(u_1+z_1,\ldots,u_d+z_d)$ for $U$, $U'$ where
$z_i\in\rad V$, and so $U\simeq U'$.
As ever, we say a nonzero vector $x\in V$
is isotropic if $Q'(x)=0$, and $V$ is isotropic if it contains a
(nonzero) isotropic
vector (and anisotropic otherwise).  We say $V$ is totally isotropic
if all its nonzero vectors are isotropic.
A hyperbolic plane is a dimension 2 space $\H$ with a basis $x,y$ so that
$Q'(x)=Q'(y)=0$, $B(x,y)=1$.

In the following we make frequent use of the fact that $\F^2=\F$
(recall $\gamma\mapsto\gamma^2$ is a homomorphism from $\F^*$ to $\F^*$
with kernel $\{1\}$).
We also make use of the consequence that an anisotropic line represents
every element of $\F$ exactly once.

\proclaim{Proposition 5.1} Say $V$ is regular with $\dim V\ge 3$.  Then $V$
is istropic.\endproclaim

\demo{Proof}  It suffices to consider $\dim V=3$.  Let $x_1, x_2, x_3$
be a basis for $V$.  If $x_i$ is isotropic for some $i$, then we are done.
So say $Q'(x_i)\not=0$ for all $i$.
 If $B(x_1,x_2)=0$ then we can choose $\gamma$
so that $\gamma x_1+ x_2$ is isotropic.
So suppose $B(x_1,x_2), B(x_1,x_3)\not=0$.  Then we can choose
$\delta$ so that $B(x_1,x_2+\delta x_3)=0$, and then we can choose
$\gamma$ so that $\gamma x_1+(x_2+\delta x_3)$ is isotropic. $\square$
\enddemo

\proclaim{Proposition 5.2} Say $V$ is regular and 
$x\in V$ is isotropic.  Then $x$ lies in a hyperbolic plane that splits $V$.
\endproclaim

\demo{Proof}  Since $V$ is regular, there is
some $y\in V$ so that $B(x,y)\not=0$.  We can scale
$y$ so that $B(x,y)=1$.  If $Q'(y)=0$ then $\F x\oplus\F y$ is a hyperbolic plane.
So suppose $Q'(y)\not=0$.  Then $Q'(y)x+y$ is isotropic, and
$\F x\oplus\F y=\F x\oplus\F (Q'(y)x+y)$ is a hyperbolic plane. 

Let $x,y,z_1,\ldots,z_d$ be a basis for $V$, $Q'(x)=Q'(y)=0$, $B(x,y)=1$.
Set $z'_i=z_i+B(z_i,y)x+B(z_i,x)y$.  Then $x,y,z'_1,\ldots,z'_d$ is also
a basis for $V$, and each $z'_i$ is orthogonal to $x$ and $y$.  
Thus $\F x\oplus\F y$ splits $V$.  $\square$
\enddemo

Note that these two propositions imply
 that for $V$ regular, $V\simeq \H^d\perp W$
where $W$ is anisotropic of dimension 0, 1, or 2.  
Since $\F^2=\F$, any two 1-dimensional regular spaces are isometric, and 
hence any two regular spaces with dimension $2d+1$ are isometric.

The following proposition
shows that an anisotropic plane, i.e. an anisotropic space of dimension 2,
cannot be diagonal.

\proclaim{Proposition 5.3}  Say $V$ is regular with an orthogonal basis
$x_1,\ldots,x_m$.  Then $m=1$.
\endproclaim

\demo{Proof}  It suffices to show that $\F x\perp\F y$ is not regular.
If either $x$ or $y$ is isotropic then this is clear.  Otherwise, we
can choose $\gamma\in\F$ so that 
$Q'(x+\gamma y)=Q'(x)+\gamma^2 Q'(y)=0$ (recall $\F^2=\F$). 
Then we also have $B(x,x+\gamma y)=B(y,x+\gamma y)=0$, showing
$V$ is not regular.
$\square$
\enddemo

In what follows we use $H=\{\gamma^2+\gamma:\ \gamma\in\F\ \}.$
Since $\gamma\mapsto\gamma^2+\gamma$ is a homomorphism of the additive
group $\F$ with kernel $\{0,1\}$, $H$ is an additive subgroup of $\F$
with index 2.

\proclaim{Proposition 5.4}  Let $W$ be a regular plane.
Then $W$ has a basis $x, y$ so that $Q'(x)=B(x,y)=1$, and
$W$ is a hyperbolic plane if and only if $Q'(y)\in 
H=\{\gamma^2+\gamma:\ \gamma\in\F\ \}.$  Further, for $W$ anisotropic
and $\varepsilon\not\in H$, $W$ has a basis $x, y$ so that
$Q'(x)=B(x,y)=1$, $Q'(y)=\varepsilon$.
\endproclaim

\demo{Proof}  Let $v, w$ be vectors so that $W=\F x\oplus\F y$.  
By Proposition 5.3,
we know $B(x,y)\not=0$.  Thus by swapping $x$ and $y$, or by replacing
$x$ by $x+y$, we can assume $Q'(x)\not=0$.  Then, since $\F^2=\F$,
we can scale $x$ to assume $Q'(x)=1$, then scale $y$ to assume $B(x,y)=1$.
Let $\varepsilon=Q'(y)$; we see $W$ is isotropic if and only if
$\alpha^2+\alpha\beta+\beta^2\varepsilon=0$ for some $\alpha,\beta$ not
both 0.  Thus $W$ is isotropic if and only if $\varepsilon=0$ or
$\varepsilon=(\alpha/\beta)^2+(\alpha/\beta)$ for some $\alpha,\beta\in\F$.
Hence by this together with Proposition 5.2, $W$ is hyperbolic if and only
$\varepsilon\in H$.

Now say $W$ is anisotropic and $\varepsilon\not\in H$.  We have
$W=\F x\oplus\F y$, $Q'(x)=B(x,y)=1$, $Q'(x)=B(x,y)=1$, $Q'(y)\not\in H$.
Since $H$ has index 2 in $\F$, we have $\varepsilon+Q'(y)\in H$,
so $\varepsilon+Q'(y)=\gamma^2+\gamma$ for some $\gamma\in\F$.  Thus
$W=\F x\oplus\F(\gamma x+y)$, with $Q'(x)=B(x,\gamma x+y)=1$, 
$$Q'(\gamma x+y)=\gamma^2 Q'(x)+\gamma B(x,y)+Q'(y)=\varepsilon.
\ \square$$
\enddemo

Note that this proposition shows all anisotropic planes are isometric,
and thus we simply use $\A$ to denote an anisotropic plane.

When working over a finite field of odd characteristic, we determine
the structure of an orthogonal sum of regular planes
$W_1\perp\cdots\perp W_m$ by determining whether $(-1)^m dW_1\cdots dW_m$
is a square; here $dW$ denotes the discriminant of $W$, or equivalently,
the determinant of a symmetric matrix representing the quadratic form
on $W$.  In the characteristic 2 case, 
Proposition 5.4 (above) and Proposition 5.5 (below)
 show that the structure of
an orthogonal sum of regular planes
$V=W_1\perp\cdots\perp W_m$ is determined by whether
$$Q'(y_1)+\cdots+Q'(y_m)\in H=\{\gamma^2+\gamma:\ \gamma\in\F\ \},$$
where $W_i=\F x_i+\F y_i$,
$Q'(x_i)=B(x_i,y_i)=1$; if $Q'(y_1)+\cdots+Q'(y_m)\in H$ then
$V\simeq \H^m$, and otherwise $V\simeq \H^{m-1}\perp\A$.

\proclaim{Proposition 5.5}  With $\A$ an anisotropic plane and $\H$
a hyperbolic plane, $\A\perp\A\simeq\H\perp\H$.
\endproclaim

\demo{Proof}  Say $W=\F x\oplus\F y$, $W'=\F x'\oplus\F y'$ are anisotropic planes;
so by Proposition 5.4
we can assume $Q'(x)=Q'(x')=B(x,y)=B(x',y')=1$, $Q'(y)=Q'(y')=\varepsilon$
for some $\varepsilon\not\in H$.  
Then $W\perp W'=U\perp U'$ where $U=\F(x+x')\oplus\F y$, 
$U'=\F x'\oplus\F(y+y')$.
Proposition 5.2 shows that $U, U'$ are hyperbolic planes.
Thus we have $$\A\perp\A\simeq W\perp W'=U\perp U'\simeq\H\perp\H.
\ \square$$
\enddemo

In Theorem 5.7 we prove a ``cancellation'' theorem, showing $\H$ can
be cancelled across an isometry.  
We first establish the following lemma.
Note that the proofs of Lemma 5.6 and Theorem 5.7 mimic the proofs of
Proposition 93:12 and Theorem 93:14 of [8].

\proclaim{Lemma 5.6}  Say $V=W\perp U$, $W=\F x \oplus \F y$ a hyperbolic
plane with $Q'(x)=Q'(y)=0$, $B(x,y)=1$.  Say $z\in U$; set
$W'=\F (x+z) \oplus \F y$.  Then $W'\simeq \H$ and $V=W'\perp U'$
with $U\simeq U'$.
\endproclaim

\demo{Proof}  By Proposition 5.2, $W'\simeq \H$.  We define
$\sigma: U\to V$ by 
$$\sigma(u)=u+B(u,z)y.$$
Thus $\sigma$ is a linear transformation, and since $u,y$ are linearly
independent when $u\not=0$, $\sigma$ is injective.
Also, recalling that $W$ is orthogonal to $U$, we find
$Q'(\sigma u)=Q'(u),$ and $B(\sigma u,x+z)=0=B(\sigma u,y).$
Thus $\sigma$ is an isometry taking $U$ into $U'$ (since $U'$ is the
orthogonal complement of $W'=\F(x+z) \oplus \F y$).  Since $U, U'$
have the same (finite) cardinality, we must have $\sigma U=U'$.
$\square$
\enddemo

\proclaim{Theorem 5.7}  Say $V$ is regular, 
and $W$ is a subspace with $W\simeq\H$,
a hyperbolic plane.  Then $V=W\perp V'$
where $V'$ is regular
with $\dim V'=\dim V-2$, and $V'$ is hyperbolic if and only if
$V$ is.  In fact, if
$V\simeq W\perp U$ and
$V\simeq W'\perp U'$, $W\simeq W'\simeq\H$, then $U\simeq U'$.
\endproclaim

\demo{Proof}   We prove the second statement; examining the structure of $U$
along the way proves the first statement.

 Choose isotropic $x, y, x', y'$ so that 
$B(x,y)=B(x',y')=1$, and $W=\F x\oplus\F y$, $W'=\F x'\oplus\F y'$.

(a)  Say $x=x'$.  We can realize $y'$ as $\gamma x+\delta y+z$,
$z\in U$; since $B(x,y')=B(x',y')=1$, we must have $\delta=1$.
So
$$W=\F x\oplus\F(\gamma x+y),\qquad W'=\F x\oplus\F(\gamma x+y+z).$$
By Lemma 5.6, $U\simeq U'$.

(b)  Say $W$ is not orthogonal to $W'$; we reduce the problem to case (a).
By assumption, $B(\gamma x+\delta y,\gamma' x'+\delta' y')\not=0$
for some $\gamma,\delta,\gamma',\delta'$; without loss of generality
we can assume $B(x,y')=1$.  Set $W''=\F x\oplus\F y'$, a hyperbolic plane.
So by Proposition 5.2, $V=W''\perp U''$, and by (a),
$U\simeq U''$ and $U'\simeq U''$.  Hence $U\simeq U'$.

(c)  Say $W$ is orthogonal to $W'$; we reduce the problem to case (b).
Set $W''=\F x\oplus\F(y+y')$, a hyperbolic plane.  Note that $W, W'$ are not
orthogonal to $W''$, since $B(x,y+y')=1=B(x',y+y')$.  So by (b),
$U\simeq U''\simeq U'$.  $\square$
\enddemo

In Theorem 5.9 we prove that one can cancel an anisotropic plane across
an isometry.  To prove this
we need the following, which is an approximation of being able to split a 
space with an anisotropic vector when $\dim V$ is even.

\proclaim{Proposition 5.8}  Say $V$ is regular and isotropic of even dimension,
and $v$
is an anisotropic vector of $V$.  Then $v$ lies in a hyperbolic plane
$U$ that splits $V$ as $U\perp V'$, and 
$v^{\perp}=\F v\perp V'$.
\endproclaim

\demo{Proof}
First  say  $V$ is hyperbolic.  So
$$V=(\F x_1\oplus\F y_1)\perp \cdots \perp (\F x_d\oplus\F y_d),$$
where $x_i,y_i$ are isotropic with $B(x_i,y_i)=1$.
Hence $v=x+y$ where $x=\sum_i \gamma_i x_i$, $y=\sum_i \delta_i y_i$
are isotropic and $B(x,y)=Q'(v)\not=0$.  
Thus $v\in \F x\oplus\F y\simeq\H$.

Now say
$$V=(\F x_1\oplus\F y_1)\perp \cdots \perp (\F x_d\oplus\F y_d)\perp W,$$
where $x_i,y_i$ are isotropic with $B(x_i,y_i)=1$, and $W\simeq\A$ 
(an anisotropic
plane).  So $v=x+y+w$ where
$x=\sum_i\gamma_i x_i$, $y=\sum_i\delta_i y_i$, $w\in W$.  If $w=0$
then this reduces to the preceeding case.  So say $w\not=0$;
hence $Q'(w)\not=0$.  Say $x\not=0$; then $\gamma_i\not=0$ for some $i$,
so $B(v,y_i)=\gamma_i\not=0$ and hence $\F v\oplus\F y_i\simeq\H$, proving
the claim.
So now suppose $x=y=0$.  Choose $u\in W$ so that $\F w\oplus\F u=W$;
since $W$ is anisotropic and $\F^2=\F$, we can scale $u$ to assume $Q'(u)=1$.
Then $x_1+y_1+u$ is isotropic, 
$B(v, x_1+y_1+u ) = B(v,u)\not=0,$
so $\F w\oplus\F(x_1+y_1+u) \simeq \H.$
$\square$
\enddemo

\proclaim{Theorem 5.9} Say $V$ is regular and $W$ is a subspace with
$W\simeq\A,$ an anisotropic plane.  Then $V=W\perp V'$
where $V'$ is regular; when $\dim V$ is even, $V'$ is
hyperbolic if and only if $V$ is not.
In fact, if $V=U\perp W=U'\perp W'$ with $W\simeq W'\simeq\A$,
then $U\simeq U'$.
\endproclaim

\demo{Proof}
We first prove the theorem for $\dim V$ even; then, using this, we prove the
theorem for $\dim V$ odd.

So we first suppose
that $\dim V=2m$, $m\ge 1$; if $m=1$ then $V=W$
and we are done.
So suppose $m>1$; thus $V$ is isotropic by Proposition 5.1.
By Proposition 5.4 we know that for some $v, w$,
 $W=\F v\oplus \F w$ where $Q'(v)=B(v,w)=1$, $Q'(w)\not\in H=
\{\gamma^2+\gamma:\ \gamma\in\F\ \}$.
So by Proposition 5.8,
$v=x+y$ where $x, y$ are isotropic,
$B(x,y)=Q'(v)=1$. Then by Theorem 5.7,
$V=(\F x\oplus\F y)\perp U$, $U$ regular of dimension $2(m-1)$, and
$U$ hyperbolic if and only if $V$ is.
Thus $w=\alpha x+\beta y+u$ for some $u\in U$ and $\alpha,\beta\in\F$; so
$\alpha+\beta=B(v,w)=1$.
Hence,  noting that $\alpha+1=\beta$, we have $w=\alpha(x+y)+(y+u)$ and so
$$W = \F(x+y) \oplus \F(y+u),$$
where $Q'(x+y)=B(x+y,y+u)=1$, $Q'(y+u)=Q'(u)$.

To prove the second statement of the theorem when $\dim V$ is even,
we separate the cases of $U$ isotropic and $U$ anisotropic; we continue
to assume $w=\alpha x+\beta y+u$ with conditions as above.

First suppose $U$ is isotropic.  Then by Proposition 5.8,
$u=x'+y'$ where $x', y'$ are isotropic with $B(x',y')=Q'(u)$.
Hence by Theorem 5.7, $V=(\F x\oplus\F y)\perp(\F x'\oplus\F y')\perp U'$
where $U'$ is regular of dimension of $2(m-2)$, and hyperbolic if and only if
$V$ is.  So $W=\F(x+y)\oplus\F(y+x'+y')$, and
$W^{\perp}=W'\perp U'$ where $W'=\F(x+y+\gamma x')\oplus\F(x'+y')$,
$\gamma Q'(u)=1$.  So $V=W\perp W'\perp U'$.  By direct computation,
or by Proposition 5.5 and Theorem 5.7, $W'\simeq \A$.

Next say $U$ is anisotropic.  Thus $m=2$ and $U$ is an anisotropic
plane.  Choose $z\in U$ so that $U=\F u\oplus\F z$ with $B(u,z)=1$.
Thus $W=\F(x+y)\oplus\F(y+u)$ and $W^{\perp}=\F(x+y+z)\oplus\F u$.
Hence $U=W\perp W^{\perp}$.  By direct computation, or by Proposition 5.5
and Theorem 5.7, $W^{\perp}\simeq\H$.  

For $\dim V$ even,
the last claim of the theorem now follows easily from Propositions 5.1,
5.2, and 5.4.

Now suppose $\dim V$ is odd.  Thus, as discussed following Proposition 5.2,
$V=V'\perp V_0$ where $V'\simeq\H^d$, $V_0$ is anisotropic, and necessarily
$\dim V_0=1$.  Since any vector (and hence any line) is orthogonal to itself,
$$V_0=\{u\in V:\ B(u,V)=0\ \}.$$
 Note that with $W\simeq\A$, $W\subseteq V$,
we must have $W\cap V_0=\{0\}$.  Thus we can decompose $V$ as
$V=V_0\perp V'$ where $W\subseteq V'$, and $V'$ is regular of even dimension.
Then by the previous part of this proof, $W$ splits $V'$, and so $W$ splits
$V$.  Also, if $\dim V$ is odd and
$V=U\perp W=U'\perp W'$ with $W\simeq W'\simeq\A$, then 
$U, U'$ are regular with $\dim U=\dim U'$
odd, and so by Propositions 5.1 and 5.2, $U\simeq U'$.
$\square$
\enddemo

We finish this section with two results regarding totally isotropic
subspaces.

\proclaim{Proposition 5.10}  Say $V$ is regular and $R$ is a totally isotropic
subspace of $V$ of dimension $r$.
  Then $V$ contains an $r$-dimensional subspace $R'$
so that $R\oplus R'\simeq\H^r$ (and thus $R\oplus R'$ splits $V$).
\endproclaim

\demo{Proof}  We argue by induction on $r$.  If $r=1$ then the result
follows by Proposition 5.2.  

So say $r>1$.  Choose isotropic $x\in R$.
Since $V$ is regular,
there is some $y\in V$ so that $B(x,y)\not=0$ and so by Proposition 5.2,
$\F x\oplus\F y\simeq\H$; hence by Theorem 5.7,
$V=(\F x\oplus\F y)\perp V'$, $V'$ regular.
Also, we can choose totally isotropic
$W\subseteq V'$ of dimension $r-1$ so that
$R=\F x\perp W$.
Induction gives us that $V'$ contains $W\oplus W'\simeq \H^{r-1}$.
Setting $R'=\F y\perp W'$ completes the proof.
$\square$
\enddemo

\noindent{\bf Remark.}  Say $V$ is regular of dimension $2\ell$, and $U$
is a subspace of dimension $d$.  Thus $U=R\perp W$, $R=\rad U$, $W$ regular
of dimension $d-r$, $r=\dim R$.  By the above proposition, there is a
dimension $r$ subspace $R'$ in $V$ so that
$R\oplus R'\simeq\H^r$ and $V=(R\oplus R')\perp V'$, $V'$ regular of dimension
$2(\ell-r)$, and $V'$ hyperbolic if and only if $V$ is.
Since $U=R\perp W$ and
$W\subseteq R^{\perp}=R\perp V'$, we can 
adjust $W$ to assume $W\subseteq V'$.  By
Propositions 5.1, 5.2, $W\simeq\H^t\perp W'$ where
$W'$ is anisotropic of dimension 0, 1 or 2.
(So $d=r+2t+\dim W'$.)

Say $\dim W'=0$.  Then by Theorem 5.7, $V'=W\perp V''$ where $V''$ is regular
of dimension $2(\ell-r-t)$, with $V''$ hyperbolic if and only if
$V$ is.  (Note that $V'$ is hyperbolic if and only if $V$ is.)  Then 
$$U^{\perp}=R\perp V''\simeq
\cases R\perp\H^{\ell-r-t}&\text{if $V$ is hyperbolic,}\\
R\perp\H^{\ell-r-t-1}\perp\A&\text{otherwise.}\endcases$$

Next, say $\dim W'=2$.  Then using Theorem 5.9 and arguing essentially
as above, we find
$$U^{\perp}\simeq\cases R\perp\H^{\ell-r-t-1}\perp\A
&\text{if $V$ is hyperbolic,}\\
R\perp\H^{\ell-r-t}&\text{otherwise.}\endcases$$

Finally, say $\dim W'=1$.  So $W'=\F x$, $x$ anisotropic.  Thus by
Proposition 5.2 and Theorem 5.7, there exists $y\in V'$ so that
$$V'\simeq\H^t\perp(\F x\oplus\F y)\perp V''$$
 where $V''$ is regular
 of dimension $2(\ell-r-t=1)$, and hyperbolic if and only if $V$ is.
Then $U^{\perp} = R\perp \F x\perp V''.$

In all cases, 
$$U^{\perp}\simeq\cases U\perp\H^{\ell-d}&\text{if $V$ is hyperbolic,}\\
U\perp\H^{\ell-d-1}\perp\A&\text{otherwise.}\endcases$$
(Recall that $U\simeq R\perp\H^t\perp W'$ where $R$ is totally isotropic
of dimension $r$ and $W'$ is anisotropic.  Thus by the preceeding results
in this section,
 the expression given above for $U^{\perp}$ is meaningful
with well-defined isometry class.)

\proclaim{Theorem 5.11} Suppose $V$ is regular of dimension $m$, and
let $\varphi_{\ell}(V)$ denote the number of 
totally isotropic $\ell$-dimensional subspaces
of $V$.  Then
$$\varphi_{\ell}(V)
=\cases \beta(t,\ell)\delta(t-1,\ell)
&\text{if $\dim V=2t$ and $V$ is hyperbolic,}\\
\beta(t-1,\ell)\delta(t,\ell)
&\text{if $\dim V=2t$ and $V$ is not hyperbolic,}\\
\beta(t,\ell)\delta(t,\ell)
&\text{if $\dim V=2t+1$.}
\endcases$$
Here $\delta(m,r)=\prod_{i=0}^{r-1}(q^{m-i}+1)$, 
 $\mu(m,r)=\prod_{i=0}^{r-1}(q^{m-i}-1)$, and
$\beta(m,r)=\mu(m,r)/\mu(r,r)$ ($m,r\ge 0$).

\endproclaim

\demo{Proof}  We first consider $\ell=1$.
We let $\psi(V)=\varphi_1(V)(q-1);$
so $\psi(V)$ is the number of isotropic vectors in $V$.  We derive the
formula for $\varphi_1(V)$ by proving the corresponding formula holds for
$\psi(V)$.  We argue by induction on $d$ where
$V\simeq\H^d \perp W$
and $W$ is anisotropic of dimension 0, 1, or 2.

For $d=0$, $\psi(V)=0$, consistent with the formula claimed.
So say $U\simeq \H^d\perp W$, $d\ge 0$, $U'\simeq\H$, $V=U\perp U'$.
Given isotropic $v\in V$, we have $v=u+u'$ where $u\in U$, $u'\in U'$.
If $Q'(u)=0$ then $Q'(u')=0$; note that since $v\not=0$ we can have
$u=0$ or $u'=0$, but not both.  So the number of isotropic $v=u+u'$ with
$Q'(u)=Q'(u')=0$ is
$$\left( \psi(U)+1\right) \left( \psi(U')+1\right) -1.$$
Say $Q'(u)\not=0$.  We know $U'$ contains $(q^2-1)/(q-1)=q+1$ lines,
two of which are isotropic.  Each anisotropic line represents every element
of $\F$ exactly once (since $\F^2=\F$), so $U'$ represents any $\gamma\not=0$
exactly $q-1$ times.  Thus the number of isotropic $v=u+u'$ with
$Q'(u)\not=0$ is
$$\left(q^{\dim U}-\psi(U)-1\right) (q-1).$$
Thus
$$\psi(V) = q\psi(U) + \left(q^{\dim V}+1\right) (q-1).$$
Induction on $d$ yields the result.

Now suppose $\ell\ge1$.
Let $\psi_{\ell}(V)$ be the number of isotropic, orthogonal $\ell$-tuples
of vectors $(x_1,\ldots,x_{\ell})$ so that $x_1,\ldots,x_{\ell}$ are
linearly independent in $V$; we use induction on $\ell$ to show that
$$\psi_{\ell}(V)
=\cases 
q^{\ell(\ell-1)/2}\mu(t,\ell) \delta(t-1,\ell)
&\text{if $\dim V=2t$ and $V$ is hyperbolic,}\\
q^{\ell(\ell-1)/2}\delta(t,\ell) \mu(t-1,\ell)
&\text{if $\dim V=2t$ and $V$ is not hyperbolic,}\\
q^{\ell(\ell-1)/2}\mu(t,\ell) \delta(t,\ell)
&\text{if $\dim V=2t+1$.}
\endcases$$
We have established this for $\ell=1$ in the preceeding prargraph,
so suppose $\ell>1$ and the formula holds for $\psi_{\ell-1}(U)$,
$U$ regular.  So choose $x_1$ isotropic in $V$; we have $\psi(V)$ choices
for $x_1$.  Since $V$ is regular, there is some $y_1\in V$ so that
$B(x_1,y_1)\not=0$, and hence 
(by Proposition 5.2 and Theorem 5.7) $\F x_1+\F y_1\simeq\H$
and $V=(\F x_1+\F y_1)\perp U$, $U$ regular of dimension $m-2$, and
$U$ hyperbolic if and only if $V$ is.
Thus with $x_1$ fixed, any $\ell$-tuple of isotropic, orthogonal, linearly
independent vectors $(x_1,\ldots,x_{\ell})$ has $x_i=\gamma_1 x_1+x_i$
for some $\gamma_i\in\F$, $x'_i\in U$ ($i\ge2$).
Therefore there are $q^{\ell-1}\psi_{\ell-1}(U)$ such $\ell$-tuples
with $x_1$ prescribed.  Hence
$$\psi_{\ell}(V) = \psi(V)\cdot q^{\ell-1}\psi_{\ell-1}(U).$$
Substituting the formula for $\psi_{\ell-1}(U)$ proves the formula
for $\psi_{\ell}(V)$.

Finally, since $\psi_{\ell}(V)$ tells us how many ways we can choose
an (ordered) basis for a totally isotropic dimension $\ell$ subspace
of $V$,
$$\varphi_{\ell}(V)=\psi_{\ell}(V)/
\big[(q^{\ell}-1)(q^{\ell}-q)\cdots(q^{\ell}-q^{\ell-1})\big]$$
since $(q^{\ell}-1)(q^{\ell}-q)\cdots(q^{\ell}-q^{\ell-1})
=q^{\ell(\ell-1)/2}\mu(\ell,\ell)$
is the number of (ordered) bases of an $\ell$-dimensional space.
The formula for $\varphi_{\ell}(V)$ now easily follows.  $\square$
\enddemo

\enddocument